\documentclass{amsart}

\usepackage{amssymb,amsmath,amsthm,latexsym,graphicx}

  \theoremstyle{definition}
  \newtheorem{definition}{Definition}
  \newtheorem{notation}[definition]{Notation}
  \newtheorem{example}[definition]{Example}
  \newtheorem{remark}[definition]{Remark}
  \newtheorem{algorithm}[definition]{Algorithm}

  \theoremstyle{plain}
  \newtheorem{lemma}[definition]{Lemma}
  \newtheorem{proposition}[definition]{Proposition}
  \newtheorem{theorem}[definition]{Theorem}
  \newtheorem{corollary}[definition]{Corollary}
  \newtheorem{conjecture}[definition]{Conjecture}

\begin{document}


\title{Special identities for quasi-Jordan algebras}

\author{Murray R. Bremner}

\address{Department of Mathematics and Statistics, University of Saskatchewan,
Canada}

\email{bremner@math.usask.ca}

\author{Luiz A. Peresi}

\address{Department of Mathematics, University of S\~ao Paulo, Brazil}

\email{peresi@ime.usp.br}

\begin{abstract}
Semispecial quasi-Jordan algebras (also called Jordan dialgebras) are defined
by the polynomial identities
  \[
  a(bc) = a(cb), \qquad (ba)a^2 = (ba^2)a, \qquad (b,a^2,c) = 2(b,a,c)a.
  \]
These identities are satisfied by the product $ab = a \dashv b + b \vdash a$ in
an associative dialgebra. We use computer algebra to show that every identity
for this product in degree $\le 7$ is a consequence of the three identities in
degree $\le 4$, but that six new identities exist in degree 8. Some but not all
of these new identities are noncommutative preimages of the Glennie identity.
\end{abstract}

\maketitle


\section{Introduction}

Loday \cite{Loday2, Loday3} introduced a new variety of algebras with two
binary operations.

\begin{definition}
An \textbf{associative dialgebra} is a vector space with bilinear operations $a
\dashv b$ and $a \vdash b$, the \textbf{left} and \textbf{right} products,
satisfying these polynomial identities:
  \allowdisplaybreaks
  \begin{alignat*}{3}
  ( a \vdash b ) \vdash c &= ( a \dashv b ) \vdash c,
  &\quad
  a \dashv ( b \dashv c ) &= a \dashv ( b \vdash c ),
  \\
  ( a \dashv b ) \dashv c &= a \dashv ( b \dashv c ),
  &\quad
  ( a \vdash b ) \vdash c &= a \vdash ( b \vdash c ),
  &\quad
  ( a \vdash b ) \dashv c &= a \vdash ( b \dashv c ).
  \end{alignat*}
Since $( a \dashv b ) \vdash c = a \dashv ( b \vdash c )$ does not hold in
general, this gives a class of algebras that are ``nearly associative''; see
Shirshov \cite{Shirshov3}, Zhevlakov et al.~\cite{Zhevlakov}.
\end{definition}

\begin{definition}
A \textbf{dialgebra monomial} on the set $X$ of generators is a product $w =
\overline{\, a_1 \cdots a_n}$ where $a_1, \hdots, a_n \in X$ and the bar
indicates some placement of parentheses and some choice of operations. We
define $c(w)$, the \textbf{center} of $w$, inductively: If $w \in X$ then $c(w)
= w$; otherwise $c( w_1 \dashv w_2 ) = c(w_1)$ and $c( w_1 \vdash w_2 ) =
c(w_2)$.
\end{definition}

\begin{lemma} \emph{(Loday \cite{Loday3}, 1.7 Theorem)} \label{normalform}
If $w = \overline{\, a_1 \cdots a_n}$ and $c(w) = a_k$ then
  \[
  w =
  ( a_1 \vdash \cdots \vdash a_{k-1} )
  \vdash a_k \dashv
  ( a_{k+1} \dashv \cdots \dashv a_n ).
  \]
\end{lemma}

\begin{definition}
The \textbf{normal form} of $w$ in Lemma \ref{normalform} will be abbreviated
as
  \[
  w = a_1 \cdots a_{k-1} \widehat a_k a_{k+1} \cdots a_n.
  \]
\end{definition}

\begin{lemma} \label{basislemma} \emph{(Loday \cite{Loday3}, 2.5 Theorem)}
The monomials $a_1 \cdots a_{k-1} \widehat a_k a_{k+1} \cdots a_n$ with $k = 1,
\dots, n$ and $a_1, \hdots, a_n \in X$ form a basis of the free associative
dialgebra on $X$.
\end{lemma}

\begin{notation}
We write $FD_n$ for the multilinear subspace of degree $n$ in the free
associative dialgebra on $n$ generators.  Lemma \ref{basislemma} implies that
$\dim FD_n = n(n!)$.
\end{notation}

\begin{definition}
(Vel\'asquez and Felipe \cite{VelasquezFelipe1}) The \textbf{quasi-Jordan
product} in a dialgebra over a field of characteristic $\ne 2$ is defined by
 \[
 a \triangleleft b = \tfrac12 ( a \dashv b + b \vdash a ).
 \]
If $D$ is a dialgebra, then its \textbf{plus algebra} $D^+$ has the same
underlying vector space but the operation $a \triangleleft b$. We omit the
symbol $\triangleleft$ and the coefficient $\tfrac12$ and write
  \[
  ab = a \dashv b + b \vdash a.
  \]
\end{definition}

\begin{definition} \label{quasijordanidentities}
We consider \textbf{right commutativity}, the \textbf{quasi-Jordan identity},
and the \textbf{associator-derivation identity}:
  \[
  a ( b c ) = a ( c b ),
  \qquad
  ( b a ) a^2 = ( b a^2 ) a,
  \qquad
  ( b, a^2, c ) = 2 (b,a,c) a,
  \]
where $(a,b,c) = (ab)c - a(bc)$. The multilinear forms of the last two
identities are
  \allowdisplaybreaks
  \begin{align*}
  J &=
  (a(bc))d + (a(bd))c + (a(cd))b - (ab)(cd) - (ac)(bd) - (ad)(bc),
  \\
  K &=
  ((ab)d)c + ((ac)d)b - (a(bc))d - (a(bd))c - (a(cd))b + a((bc)d).
  \end{align*}
\end{definition}

\begin{remark}
Equivalent identities for the opposite product appear in the work of
Kolesnikov. If we replace $xy$ by $yx$ in equations (26) and (27) of
\cite{Kolesnikov}, replace $x_1,x_2,x_3,x_4$ by $a,b,c,d$ and apply right
commutativity, then we obtain
  \begin{align*}
  L &= ((ac)b)d + ((ad)b)c - (ab)(cd) - (ac)(bd) - (ad)(bc) + a((cd)b),
  \\
  M &= (b(cd))a + (b(ac))d + (b(ad))c - (ba)(cd) - (bd)(ac) - (bc)(ad).
  \end{align*}
One easily verifies that
  \begin{align*}
  &
  J(a,b,c,d) = M(b,a,d,c),
  \qquad
  K(a,b,c,d) = L(a,d,b,c) - M(b,a,c,d),
  \\
  &
  L(a,b,c,d) = J(a,b,c,d) + K(a,c,d,b),
  \qquad
  M(a,b,c,d) = J(b,a,c,d).
  \end{align*}
See also the remarks on Jordan dialgebras in Pozhidaev \cite{Pozhidaev},
Section 3.
\end{remark}

\begin{lemma}
\emph{(Vel\'asquez and Felipe \cite{VelasquezFelipe1})} The quasi-Jordan
product in an associative dialgebra satisfies right commutativity and the
quasi-Jordan identity.
\end{lemma}

\begin{lemma}
\emph{(Bremner \cite{Bremner})} The quasi-Jordan product in an associative
dialgebra satisfies the associator-derivation identity.  The identities of
Definition \ref{quasijordanidentities} imply every identity of degree $\le 4$
for the quasi-Jordan product in an associative dialgebra.
\end{lemma}

\begin{definition}
A \textbf{quasi-Jordan algebra} is a nonassociative algebra over a field of
characteristic $\ne 2, 3$ satisfying right commutativity and the quasi-Jordan
identity. A \textbf{semispecial} quasi-Jordan algebra (also called a
\textbf{Jordan dialgebra}) is a quasi-Jordan algebra satisfying the
associator-derivation identity.
\end{definition}

\begin{remark} \label{jordandialgebra}
Strictly speaking, a Jordan dialgebra as defined by Kolesnikov
\cite{Kolesnikov} and Pozhidaev \cite{Pozhidaev} has two binary operations
$\dashv$ and $\vdash$ which are in fact opposite as a result of the dialgebra
version of commutativity, $x \dashv y = y \vdash x$. We simplify the notation
by using only one operation and writing this operation as juxtaposition.
\end{remark}

\begin{definition}
A quasi-Jordan algebra is \textbf{special} if it is isomorphic to a subalgebra
of $D^+$ for some associative dialgebra $D$. Every special algebra is
semispecial.
\end{definition}

Glennie \cite{Glennie1, Glennie2, Glennie3} (see also Hentzel \cite{Hentzel3})
discovered an identity satisfied by special Jordan algebras that is not
satisfied by all Jordan algebras.  In this paper we consider the analogous
question for quasi-Jordan algebras.  We use computer algebra to show that the
identities in Definition \ref{quasijordanidentities} imply every identity of
degree $\le 7$ for the quasi-Jordan product in an associative dialgebra. We
demonstrate the existence of identities in degree 8 which do not follow from
the identities of Definition \ref{quasijordanidentities}. Some but not all of
these new identities are noncommutative preimages of the Glennie identity.
These new identities are special identities in the following sense.

\begin{definition}
A \textbf{special identity} is a polynomial identity satisfied by all special
quasi-Jordan algebras but not satisfied by all semispecial quasi-Jordan
algebras.
\end{definition}

We present an explicit special identity which has three variables and is linear
in one variable; this cannot happen for special Jordan algebras by a theorem of
Macdonald \cite{Macdonald}. Our methods depend on computational linear algebra
with large matrices over a finite field, together with the representation
theory of the symmetric group. Our computations were done with C
\cite{Kernighan}, Maple \cite{Maple} and Albert \cite{Jacobs}.


\section{Preliminaries on free nonassociative algebras}

\subsection{Free right-commutative algebras}

The simplest identity satisfied by the quasi-Jordan product is right
commutativity. Our computations depend on basic facts about free
right-commutative algebras. Kurosh \cite{Kurosh} proved that every subalgebra
of an (absolutely) free nonassociative algebra is also free; Shirshov
\cite{Shirshov1} proved the same result for free commutative and free
anticommutative algebras. The referee pointed out a simple proof that this does
not hold for a free right commutative algebra $R$: the commutator $[R,R]$ is a
subalgebra with trivial multiplication (all products are zero) and hence is
certainly not free.

\begin{lemma} \label{rightcommutativelemma}
Let $w = \overline{\, a_1 \cdots a_n}$ be a nonassociative monomial where $a_1,
\hdots, a_n \in X$ and the bar denotes some placement of parentheses.
Right-commutativity implies that in any submonomial $x = yz$ we may assume
commutativity for $z$.
\end{lemma}

\begin{proof}
By induction on $n$; for $n \le 3$ the claim is immediate.  The monomial $w$
has the unique factorization $w = uv$; the inductive hypothesis implies the
claim for $u$ and $v$.  Any right factor of a submonomial of $w$ is either $v$,
a right factor of a submonomial of $u$, or a right factor of a submonomial of
$v$.  It therefore suffices to show that we may assume commutativity for $v$.
We have the unique factorization $v = xy$, and we may assume commutativity for
$y$.  Right-commutativity implies $uv = u(xy) = u(yx)$, and by induction we may
assume commutativity for $x$.
\end{proof}

Lemma \ref{rightcommutativelemma} implies the following algorithm for
generating a complete minimal set of right-commutative association types up to
a given degree $n$.

\begin{algorithm} \label{typesalgorithm}
Assume that the right-commutative association types have been generated for
degrees $\le n{-}1$. Any right-commutative type in degree $n$ has the form $w =
uv$ where $u$ is a right-commutative type in degree $n{-}i$ and $v$ is a
commutative type in degree $i$.  This induces a total order on the association
types.
\end{algorithm}

Algorithm \ref{typesalgorithm} implies the following recursive formula for the
number $R_n$ of right-commutative types in degree $n$ which requires the number
$C_n$ of commutative types.

\begin{lemma} \label{typeslemma}
We have $C_1 = R_1 = 1$, and for $n \ge 2$ we have
  \begin{alignat*}{2}
  C_n &=
  \sum_{i=1}^{\lfloor(n-1)/2\rfloor} C_{n-i} C_i + \binom{C_{n/2}{+}1}{2},
  &\qquad
  R_n &=
  \sum_{i=1}^{n-1} R_{n-i} C_i.
  \end{alignat*}
(The binomial coefficient only appears for $n$ even.)
\end{lemma}

\begin{example}
The following table gives the numbers $C_n$ and $R_n$ for $1 \le n \le 12$,
together with the Catalan number $K_n$ of all nonassociative types in degree
$n$:
  \begin{center}
  \begin{tabular}{rrrrrrrrrrrrr}
  $n$   &
  1 & 2 & 3 & 4 & 5 & 6 & 7 & 8 & 9 & 10 & 11 & 12 \\
  $C_n$ &
  1 & 1 & 1 & 2 & 3 & 6 & 11 & 23 & 46 & 98 & 207 & 451 \\
  $R_n$ &
  1 & 1 & 2 & 4 & 9 & 20 & 46 & 106 & 248 & 582 & 1376 & 3264\\
  $K_n$ &
  1 & 1 & 2 & 5 & 14 & 42 & 132 & 429 & 1430 & 4862 & 16796 &
  58786
  \end{tabular}
  \end{center}
\end{example}

\begin{lemma}
The generating functions of $C_n$ and $R_n$ are related by the equations
  \[
  C(x) = \sum_{n=1}^\infty C_n x^n,
  \qquad
  \sum_{n=1}^\infty R_n x^n = \frac{x}{1-C(x)}.
  \]
\end{lemma}

\begin{proof}
Sloane \cite{Sloane} (sequences A001190 and A085748).
\end{proof}

\begin{definition}
In the \textbf{basic monomial} for an association type in degree $n$ the
variables are the first $n$ letters of the alphabet in lexicographical order.
\end{definition}

\begin{example}
For $n = 1$ (resp.~$n = 2$) we have the single type $a$ (resp.~$ab$). For $3
\le n \le 5$ we present the basic commutative and right-commutative monomials:
  \[
  \begin{array}{lll}
  n
  &
  \text{commutative}
  &
  \text{right-commutative}
  \\
  3
  &
  (ab)c
  &
  (ab)c, a(bc)
  \\
  4
  &
  ((ab)c)d, (ab)(cd)
  &
  ((ab)c)d, (a(bc))d, (ab)(cd), a((bc)d)
  \\
  5
  &
  (((ab)c)d)e, ((ab)(cd))e, ((ab)c)(de)
  &
  (((ab)c)d)e, ((a(bc))d)e, ((ab)(cd))e,
  \\
  &
  &(a((bc)d))e, ((ab)c)(de), (a(bc))(de),
  \\
  &
  &(ab)((cd)e), a(((bc)d)e), a((bc)(de))
  \end{array}
  \]
\end{example}

\subsection{Multilinear right-commutative monomials}

Throughout most of this paper we consider only multilinear identities: in
degree $n$, the variables in each monomial are a permutation of the first $n$
letters of the alphabet.  To obtain a basis for the space of multilinear
right-commutative polynomials in degree $n$, we need a straightening algorithm
which replaces each monomial $w$ by the first monomial (in lex order) in its
equivalence class $[w]$: the set of all monomials which are equal to $w$ as a
consequence of right-commutativity. To straighten a right-commutative monomial,
it suffices to determine the symmetries of its association type.

\begin{definition}
Let $v$ be the basic monomial for a right-commutative association type. Suppose
that $v = \cdots (xy) \cdots$ contains the submonomial $xy$ where $x$ and $y$
are submonomials with the same degree and association type.  Let $w = \cdots
(yx) \cdots$ be the monomial obtained from $v$ by transposing $x$ and $y$. If
right-commutativity implies $v = w$ then this identity will be called a
\textbf{symmetry} of the association type.
\end{definition}

\begin{lemma} \label{symmetrylemma}
If a right-commutative association type in degree $n$ has $s$ symmetries, then
the number of multilinear monomials with this association type is $n!/2^s$.
\end{lemma}

\begin{proof}
Each symmetry reduces the number of monomials by a factor of 2.
\end{proof}

\begin{example}
The symmetries of the right-commutative types in degree 5:
  \[
  \begin{array}{l}
  \text{type 1:} \quad (((ab)c)d)e \; \text{has no symmetries}
  \\
  \text{type 2:} \quad ((a(bc))d)e = ((a(cb))d)e
  \qquad
  \text{type 3:} \quad ((ab)(cd))e = ((ab)(dc))e
  \\
  \text{type 4:} \quad (a((bc)d))e = (a((cb)d))e
  \qquad
  \text{type 5:} \quad ((ab)c)(de) = ((ab)c)(ed)
  \\
  \text{type 6:} \quad (a(bc))(de) = (a(cb))(de) = (a(bc))(ed)
  \\
  \text{type 7:} \quad (ab)((cd)e) = (ab)((dc)e)
  \qquad
  \text{type 8:} \quad a(((bc)d)e) = a(((cb)d)e)
  \\
  \text{type 9:} \quad a((bc)(de)) = a((cb)(de)) = a((bc)(ed)) = a((de)(bc))
  \end{array}
  \]
These types have (respectively) 0, 1, 1, 1, 1, 2, 1, 1, 3 symmetries, and
contain 120, 60, 60, 60, 60, 30, 60, 60, 15 distinct multilinear monomials, for
a total of 525.
\end{example}

\begin{notation}
We write $FRC_n$ for the multilinear subspace of degree $n$ in the free
right-commutative algebra on $n$ generators. An ordered basis of $FRC_n$
consists of the distinct right-commutative monomials in degree $n$, ordered
first by association type and then by lex order of the underlying permutation.
\end{notation}

\begin{lemma}
If $s(i)$ is the number of symmetries in association type $i$ then
  \[
  \dim FRC_n = \sum_{i=1}^{R_n} \frac{n!}{2^{s(i)}}.
  \]
\end{lemma}

\begin{proof}
This follows directly from Lemma \ref{symmetrylemma}.
\end{proof}

\begin{algorithm}
This algorithm to find the symmetries of a right-commutative association type
(represented by a basic monomial) uses a global variable \texttt{symmetrylist},
initially empty. On input $w = uv$, the primary procedure \texttt{findsymmetry}
calls itself on input $u$ and then calls the secondary procedure
\texttt{findcommutativesymmetry} on input $v$. Writing $v = xy$, the secondary
procedure calls itself on input $x$ and then on input $y$; it then checks to
see if $x$ and $y$ have the same association type, and if so it appends the
symmetry $u(xy) = u(yx)$ to \texttt{symmetrylist}. Both procedures do nothing
if the input has degree 1; this is the basis of the recursion.
\end{algorithm}

\begin{example}
The dimension $\dim FRC_n$ (the total number of multilinear right-commutative
monomials) for $1 \le n \le 9$:
  \begin{center}
  \begin{tabular}{rrrrrrrrrr}
  $n$ &\;
  1 &\; 2 &\; 3 &\; 4 &\; 5 &\; 6 &\; 7 &\; 8 &\; 9
  \\
  $\dim FRC_n$ &\;
  1 & 2 & 9 & 60 & 525 & 5670 & 72765 & 1081080 & 18243225
  \end{tabular}
  \end{center}
These values satisfy the following formula from Sloane \cite{Sloane} (sequence
A001193).
\end{example}

\begin{conjecture}
For all $n \ge 1$ we have
  \[
  \dim FRC_n \stackrel{?}{=} \frac{n(2n{-}2)!}{2^{n-1} (n{-}1)!}.
  \]
\end{conjecture}

\subsection{The expansion map and the expansion matrix}

Recall that $FRC_n$ and $FD_n$ are the spaces of multilinear right-commutative
and dialgebra polynomials.

\begin{definition}  \label{expansionmap}
We define the \textbf{expansion map} $E_n\colon FRC_n \to FD_n$ on basis
monomials and extend linearly: if $\deg w = 1$ then $E_1(w) = w$; if $w = uv$
where $\deg u = n{-}i$ and $\deg v = i$ then
  \[
  E_n(w) = E_{n-i}(u) \dashv E_i(v) + E_i(v) \vdash E_{n-i}(u).
  \]
This map is well-defined since the quasi-Jordan product is right-commutative.
\end{definition}

\begin{lemma}
The multilinear polynomial identities in degree $n$ satisfied by the
quasi-Jordan product are precisely the (nonzero) elements of the kernel of
$E_n$.
\end{lemma}

Some of these kernel identities may be consequences of identities of lower
degree. We need to distinguish the ``old'' from the ``new'' identities.

\begin{definition}
With respect to the ordered bases of $FRC_n$ and $FD_n$, we represent $E_n$ by
the \textbf{expansion matrix} $[E_n]$: we have $[E_n]_{ij} = 1$ if dialgebra
monomial $i$ occurs in the expansion of right-commutative monomial $j$, and
$[E_n]_{ij} = 0$ otherwise.
\end{definition}

The sizes of the matrices $[E_n]$ grow very rapidly. We can use
\texttt{LinearAlgebra} in Maple to compute a basis for the nullspace of $[E_n]$
over $\mathbb{Q}$ for $n \le 5$; we can use \texttt{LinearAlgebra[Modular]} to
compute a basis over $\mathbb{F}_p$ for $n \le 6$.  For $n \ge 7$ we must make
the matrices smaller; for this we use the representation theory of the
symmetric group as described in Section \ref{sectionrepresentationtheory}.

Definition \ref{expansionmap} gives a recursive algorithm for computing the
expansion of a right-commutative monomial.  To initialize $[E_n]$, we let the
column index $j$ go from left to right, compute the expansion of the
corresponding right-commutative monomial, obtain $2^{n-1}$ dialgebra monomials,
convert each dialgebra monomial to normal form and determine its row index $i$,
and set the $(i,j)$ entry of the matrix to 1.

\subsection{Lifting multilinear identities}

Let $I(x_1,\hdots,x_n)$ be a multilinear polynomial identity in degree $n$; we
want to find all its consequences in degree $n{+}1$.

\begin{definition} \label{tidealdefinition}
In the free algebra the $T$-\textbf{ideal} generated by a polynomial identity
$I$ is the smallest ideal which contains $I$ and is sent to itself by all
endomorphisms. For the endomorphism condition we introduce a new variable
$x_{n+1}$ and substitute $x_i x_{n+1}$ for $x_i$. For the ideal condition we
multiply on the left or right by $x_{n+1}$.
\end{definition}

\begin{lemma} \label{liftinglemma}
If $I(x_1,\hdots,x_n)$ is a multilinear polynomial identity in degree $n$, then
every consequence of $I$ in degree $n{+}1$ is a linear combination of
permutations of the following $n{+}2$ multilinear identities in degree $n{+}1$:
  \begin{align*}
  &
  I( x_1 x_{n+1}, x_2, \hdots, x_n ),
  \quad \hdots, \quad
  I( x_1, \hdots, x_{n-1}, x_n x_{n+1} ).
  \\
  &
  I( x_1, \hdots, x_n ) x_{n+1},
  \quad
  x_{n+1} I( x_1, \hdots, x_n ).
  \end{align*}
\end{lemma}

\begin{definition}
The identities of Lemma \ref{liftinglemma} are the \textbf{liftings} of $I$ to
degree $n{+}1$.
\end{definition}

An identity $I$ in degree $n$ produces $(n{+}2) \cdots (n{+}k{+}1)$ liftings in
degree $n{+}k$. In general, a subset of these liftings generates all the
consequences of $I$ in degree $n{+}k$. For example, the symmetries of the
right-commutative association types in degree $n$ are the liftings of
right-commutativity from degree 3. Our choice of association types eliminates
most of the consequences; only the symmetries remain. In this paper, the most
important examples are the liftings of the multilinear identities $J$ and $K$
of Definition \ref{quasijordanidentities} from degree 4 to degree $n$.


\section{Nonexistence of new identities in degree 5}

In this section we provide detailed examples of our methods; for higher degrees
the objects we work with --- polynomial identities and expansion matrices ---
become so large that it is impossible to include all the details of the
computations.

  \begin{table}
  \small
  \begin{align*}
  J(ae,b,c,d)
  &=
  ((ae)(bc))d
  {+} ((ae)(bd))c
  {+} ((ae)(cd))b
  {-} ((ae)b)(cd)
  {-} ((ae)c)(bd)
  {-} ((ae)d)(bc)
  \\
  J(a,be,c,d)
  &=
  (a((be)c))d
  {+} (a((be)d))c
  {+} (a(cd))(be)
  {-} (a(be))(cd)
  {-} (ac)((be)d)
  {-} (ad)((be)c)
  \\
  J(a,b,ce,d)
  &=
  (a(b(ce)))d
  {+} (a(bd))(ce)
  {+} (a((ce)d))b
  {-} (ab)((ce)d)
  {-} (a(ce))(bd)
  {-} (ad)(b(ce))
  \\
  &=
  (a((ce)b))d
  {+} (a(bd))(ce)
  {+} (a((ce)d))b
  {-} (ab)((ce)d)
  {-} (a(ce))(bd)
  {-} (ad)((ce)b)
  \\
  J(a,b,c,de)
  &=
  (a(bc))(de)
  {+} (a(b(de)))c
  {+} (a(c(de)))b
  {-} (ab)(c(de))
  {-} (ac)(b(de))
  {-} (a(de))(bc)
  \\
  &=
  (a(bc))(de)
  {+} (a((de)b))c
  {+} (a((de)c))b
  {-} (ab)((de)c)
  {-} (ac)((de)b)
  {-} (a(de))(bc)
  \\
  J(a,b,c,d)e
  &=
  ((a(bc))d)e
  {+} ((a(bd))c)e
  {+} ((a(cd))b)e
  {-} ((ab)(cd))e
  {-} ((ac)(bd))e
  {-} ((ad)(bc))e
  \\
  eJ(a,b,c,d)
  &=
  e((a(bc))d)
  {+} e((a(bd))c)
  {+} e((a(cd))b)
  {-} e((ab)(cd))
  {-} e((ac)(bd))
  {-} e((ad)(bc))
  \\
  &=
  e(((bc)a)d)
  {+} e(((bd)a)c)
  {+} e(((cd)a)b)
  {-} e((ab)(cd))
  {-} e((ac)(bd))
  {-} e((ad)(bc))
  \\
  K(ae,b,c,d)
  &=
  (((ae)b)d)c
  {+} (((ae)c)d)b
  {-} ((ae)(bc))d
  {-} ((ae)(bd))c
  {-} ((ae)(cd))b
  {+} (ae)((bc)d)
  \\
  K(a,be,c,d)
  &=
  ((a(be))d)c
  {+} ((ac)d)(be)
  {-} (a((be)c))d
  {-} (a((be)d))c
  {-} (a(cd))(be)
  {+} a(((be)c)d)
  \\
  K(a,b,ce,d)
  &=
  ((ab)d)(ce)
  {+} ((a(ce))d)b
  {-} (a(b(ce)))d
  {-} (a(bd))(ce)
  {-} (a((ce)d))b
  {+} a((b(ce))d)
  \\
  &=
  ((ab)d)(ce)
  {+} ((a(ce))d)b
  {-} (a((ce)b))d
  {-} (a(bd))(ce)
  {-} (a((ce)d))b
  {+} a(((ce)b)d)
  \\
  K(a,b,c,de)
  &=
  ((ab)(de))c
  {+} ((ac)(de))b
  {-} (a(bc))(de)
  {-} (a(b(de)))c
  {-} (a(c(de)))b
  {+} a((bc)(de))
  \\
  &=
  ((ab)(de))c
  {+} ((ac)(de))b
  {-} (a(bc))(de)
  {-} (a((de)b))c
  {-} (a((de)c))b
  {+} a((bc)(de))
  \\
  K(a,b,c,d)e
  &=
  (((ab)d)c)e
  {+} (((ac)d)b)e
  {-} ((a(bc))d)e
  {-} ((a(bd))c)e
  {-} ((a(cd))b)e
  {+} (a((bc)d))e
  \\
  eK(a,b,c,d)
  &=
  e(((ab)d)c)
  {+} e(((ac)d)b)
  {-} e((a(bc))d)
  {-} e((a(bd))c)
  {-} e((a(cd))b)
  {+} e(a((bc)d))
  \\
  &=
  e(((ab)d)c)
  {+} e(((ac)d)b)
  {-} e(((bc)a)d)
  {-} e(((bd)a)c)
  {-} e(((cd)a)b)
  {+} e(((bc)d)a)
  \end{align*}
  \caption{Liftings of $J$ and $K$ to degree 5}
  \label{JKliftingtable}
  \end{table}

\subsection{Old identities}

Identities $J$ and $K$ each have six liftings to degree 5. The terms of each
lifting must be straightened to lie in the standard basis of $FRC_5$; see Table
\ref{JKliftingtable}. We allocate memory for a matrix $M$ of size $645 \times
525$ with a $525 \times 525$ upper block and a $120 \times 525$ lower block;
525 is the number of multilinear right-commutative monomials and 120 is the
number of permutations of 5 variables. For each of the 12 lifted and
straightened identities $L$ in Table \ref{JKliftingtable}, we do the following:
for permutations $\pi_1, \dots, \pi_{120}$ of $a,b,c,d,e$ in lex order we apply
$\pi_j$ to $L$, straighten the terms, and store the resulting coefficient
vector in row $525{+}j$ of $M$; we then compute the row canonical form of $M$
and record the rank (the lower block of $M$ is now zero). Using rational
arithmetic we obtain the following ranks: 20, 50, 50, 50, 70, 90, 150, 210,
210, 220, 250, 250. The lifted identities which do not increase the rank are
redundant, so we consider only numbers 1, 2, 5, 6, 7, 8, 10, 11. Using modular
arithmetic we obtain the same ranks approximately 1000 times faster.

\begin{lemma}
The polynomial identities in degree 5 which are consequences of identities in
degree $\le 4$ span a 250-dimensional subspace of $FRC_5$.
\end{lemma}

  \begin{table}
  \small
  \begin{align*}
  (((ab)c)d)e
  &\longmapsto
  \widehat{a}bcde +
  e\widehat{a}bcd +
  d\widehat{a}bce +
  ed\widehat{a}bc +
  c\widehat{a}bde +
  ec\widehat{a}bd +
  dc\widehat{a}be +
  edc\widehat{a}b
  \\
  &\qquad +
  b\widehat{a}cde +
  eb\widehat{a}cd +
  db\widehat{a}ce +
  edb\widehat{a}c +
  cb\widehat{a}de +
  ecb\widehat{a}d +
  dcb\widehat{a}e +
  edcb\widehat{a}
  \\
  ((a(bc))d)e
  &\longmapsto
  \widehat{a}bcde +
  e\widehat{a}bcd +
  d\widehat{a}bce +
  ed\widehat{a}bc +
  bc\widehat{a}de +
  ebc\widehat{a}d +
  dbc\widehat{a}e +
  edbc\widehat{a}
  \\
  &\qquad +
  \widehat{a}cbde +
  e\widehat{a}cbd +
  d\widehat{a}cbe +
  ed\widehat{a}cb +
  cb\widehat{a}de +
  ecb\widehat{a}d +
  dcb\widehat{a}e +
  edcb\widehat{a}
  \\
  ((ab)(cd))e
  &\longmapsto
  \widehat{a}bcde +
  e\widehat{a}bcd +
  cd\widehat{a}be +
  ecd\widehat{a}b +
  \widehat{a}bdce +
  e\widehat{a}bdc +
  dc\widehat{a}be +
  edc\widehat{a}b
  \\
  &\qquad +
  b\widehat{a}cde +
  eb\widehat{a}cd +
  cdb\widehat{a}e +
  ecdb\widehat{a} +
  b\widehat{a}dce +
  eb\widehat{a}dc +
  dcb\widehat{a}e +
  edcb\widehat{a}
  \\
  (a((bc)d))e
  &\longmapsto
  \widehat{a}bcde +
  e\widehat{a}bcd +
  bcd\widehat{a}e +
  ebcd\widehat{a} +
  \widehat{a}dbce +
  e\widehat{a}dbc +
  dbc\widehat{a}e +
  edbc\widehat{a}
  \\
  &\qquad +
  \widehat{a}cbde +
  e\widehat{a}cbd +
  cbd\widehat{a}e +
  ecbd\widehat{a} +
  \widehat{a}dcbe +
  e\widehat{a}dcb +
  dcb\widehat{a}e +
  edcb\widehat{a}
  \\
  ((ab)c)(de)
  &\longmapsto
  \widehat{a}bcde +
  de\widehat{a}bc +
  \widehat{a}bced +
  ed\widehat{a}bc +
  c\widehat{a}bde +
  dec\widehat{a}b +
  c\widehat{a}bed +
  edc\widehat{a}b
  \\
  &\qquad +
  b\widehat{a}cde +
  deb\widehat{a}c +
  b\widehat{a}ced +
  edb\widehat{a}c +
  cb\widehat{a}de +
  decb\widehat{a} +
  cb\widehat{a}ed +
  edcb\widehat{a}
  \\
  (a(bc))(de)
  &\longmapsto
  \widehat{a}bcde +
  de\widehat{a}bc +
  \widehat{a}bced +
  ed\widehat{a}bc +
  bc\widehat{a}de +
  debc\widehat{a} +
  bc\widehat{a}ed +
  edbc\widehat{a}
  \\
  &\qquad +
  \widehat{a}cbde +
  de\widehat{a}cb +
  \widehat{a}cbed +
  ed\widehat{a}cb +
  cb\widehat{a}de +
  decb\widehat{a} +
  cb\widehat{a}ed +
  edcb\widehat{a}
  \\
  (ab)((cd)e)
  &\longmapsto
  \widehat{a}bcde +
  cde\widehat{a}b +
  \widehat{a}becd +
  ecd\widehat{a}b +
  \widehat{a}bdce +
  dce\widehat{a}b +
  \widehat{a}bedc +
  edc\widehat{a}b
  \\
  &\qquad +
  b\widehat{a}cde +
  cdeb\widehat{a} +
  b\widehat{a}ecd +
  ecdb\widehat{a} +
  b\widehat{a}dce +
  dceb\widehat{a} +
  b\widehat{a}edc +
  edcb\widehat{a}
  \\
  a(((bc)d)e)
  &\longmapsto
  \widehat{a}bcde +
  bcde\widehat{a} +
  \widehat{a}ebcd +
  ebcd\widehat{a} +
  \widehat{a}dbce +
  dbce\widehat{a} +
  \widehat{a}edbc +
  edbc\widehat{a}
  \\
  &\qquad +
  \widehat{a}cbde +
  cbde\widehat{a} +
  \widehat{a}ecbd +
  ecbd\widehat{a} +
  \widehat{a}dcbe +
  dcbe\widehat{a} +
  \widehat{a}edcb +
  edcb\widehat{a}
  \\
  a((bc)(de))
  &\longmapsto
  \widehat{a}bcde +
  bcde\widehat{a} +
  \widehat{a}debc +
  debc\widehat{a} +
  \widehat{a}bced +
  bced\widehat{a} +
  \widehat{a}edbc +
  edbc\widehat{a}
  \\
  &\qquad +
  \widehat{a}cbde +
  cbde\widehat{a} +
  \widehat{a}decb +
  decb\widehat{a} +
  \widehat{a}cbed +
  cbed\widehat{a} +
  \widehat{a}edcb +
  edcb\widehat{a}
  \end{align*}
  \caption{Expansions of the basic monomials in degree 5}
  \label{expansiontable}
  \end{table}

\subsection{All identities}

We allocate memory for the expansion matrix $E = [E_5]$ of size $600 \times
525$.  We compute the expansions of the basic monomials for the
right-commutative association types; see Table \ref{expansiontable}. From these
we obtain the expansions of all 525 monomials corresponding to the columns of
$E$; we set to 1 the appropriate entries of $E$. We obtain a sparse 0-1 matrix
in which each column has 16 nonzero entries. We find that the rank of this
matrix is 275, and so the nullity is 250.

\begin{lemma}
The subspace of $FRC_5$ consisting of polynomial identities satisfied by the
quasi-Jordan product has dimension 250.
\end{lemma}

\begin{proposition}
Every polynomial identity in degree 5 for the quasi-Jordan product is a
consequence of identities in degree $\le 4$.
\end{proposition}

\begin{proof}
The subspace generated by the lifted identities is contained in the subspace of
all identities; since the dimensions are equal, the subspaces are equal.
\end{proof}


\section{Nonexistence of new identities in degree 6: first computation}

\begin{lemma}
In degree 6 there are 20 right-commutative association types:
  \[
  \begin{array}{lllll}
  ((((ab)c)d)e)f
  & (((a(bc))d)e)f
  & (((ab)(cd))e)f
  & ((a((bc)d))e)f
  & (((ab)c)(de))f
  \\
  ((a(bc))(de))f
  & ((ab)((cd)e))f
  & (a(((bc)d)e))f
  & (a((bc)(de)))f
  & (((ab)c)d)(ef)
  \\
  ((a(bc))d)(ef)
  & ((ab)(cd))(ef)
  & (a((bc)d))(ef)
  & ((ab)c)((de)f)
  & (a(bc))((de)f)
  \\
  (ab)(((cd)e)f)
  & (ab)((cd)(ef))
  & a((((bc)d)e)f)
  & a(((bc)(de))f)
  & a(((bc)d)(ef))
  \end{array}
  \]
Each type has (respectively) 720, 360, 360, 360, 360, 180, 360, 360, 90, 360,
180, 180, 180, 360, 180, 360,  90, 360, 90, 180 monomials, for a total of 5670.
\end{lemma}

\begin{proof}
This follows directly from Lemmas \ref{typeslemma} and \ref{symmetrylemma}.
\end{proof}

For a matrix with 5670 columns, it is not practical to use rational arithmetic
to compute the row canonical form.  Instead we use modular arithmetic (with $p
= 101$) to compute the dimensions of the subspaces of lifted identities and all
identities.

\subsection{Old identities}

Our computations in degree 5 showed that we needed only 8 liftings to generate
all consequences of $J$ and $K$ in degree 5:
  \allowdisplaybreaks
  \begin{alignat*}{4}
  &J(ae,b,c,d), &\qquad &J(a,be,c,d), &\qquad &J(a,b,c,d)e, &\qquad &eJ(a,b,c,d),
  \\
  &K(ae,b,c,d), &\qquad &K(a,be,c,d), &\qquad &K(a,b,c,de), &\qquad &K(a,b,c,d)e.
  \end{alignat*}
Each of these identities produces 7 liftings in degree 6, and so we obtain an
ordered list of 56 liftings in degree 6. We follow the same algorithm as for
degree 5, except that now the matrix $M$ has size $6390 \times 5670$ with a
$5670 \times 5670$ upper block and a $720 \times 5670$ lower block. To each of
the 56 liftings, we apply all 720 permutations of the 6 variables and
straighten the terms to obtain monomials in the standard basis of $FRC_6$; we
store the coefficient vectors of the permuted liftings in the lower block, and
compute the row canonical form. We obtain the following ranks: 120, 300, 300,
300, 360, 480, 540, 540, 720, 810, 810, 810, 990, 1170, 1170, 1170, 1170, 1170,
1230, 1350, 1410, 1410, 1410, 1410, 1410, 1530, 1626, 1626, 1986, 2346, 2346,
2406, 2586, 2766, 2766, 2766, 3126, 3210, 3330, 3330, 3510, 3510, 3510, 3510,
3510, 3510, 3510, 3570, 3570, 3570, 3570, 3570, 3570, 3570, 3690, 3690. Only 25
liftings in the ordered list produce an increase in the rank: numbers 1, 2, 5,
6, 7, 9, 10, 13, 14, 19, 20, 21, 26, 27, 29, 30, 32, 33, 34, 37, 38, 39, 41,
48, 55.

\begin{lemma}
The polynomial identities in degree 6 which are consequences of identities in
degree $\le 5$ span a 3690-dimensional subspace of $FRC_6$.
\end{lemma}

\subsection{All identities}

The expansion matrix $E = [E_6]$ has size $4320 \times 5670$.  As for degree 5,
we compute the expansions of the basic monomials, determine the normal forms of
the dialgebra monomials, obtain the expansions of all the multilinear
right-commutative monomials, and store the results in the columns of $E$.  We
obtain a very sparse 0-1 matrix in which each column has 32 nonzero entries. We
find that the rank of this matrix is 1980, and so the nullity is 3690.

\begin{lemma}
The subspace of $FRC_6$ consisting of polynomial identities satisfied by the
quasi-Jordan product has dimension 3690.
\end{lemma}

\begin{proposition}
Every polynomial identity in degree 6 for the quasi-Jordan product (over the
field $\mathbb{F}_{101}$) is a consequence of identities in degree $\le 5$.
\end{proposition}

We next show how to obtain the same results using much smaller matrices.


\section{Preliminaries on representation theory}
\label{sectionrepresentationtheory}

\subsection{Representations of semisimple algebras}

Let $A$ be a finite-dimensional semisimple associative algebra over a field
$F$; then $A$ is the direct sum of simple two-sided ideals which are orthogonal
as subalgebras:
  \begin{equation} \label{semisimplealgebra}
  A = A_1 \oplus \cdots \oplus A_r,
  \qquad
  A_i A_j = \{0\} \; (1 \le i \ne j \le r).
  \end{equation}
Each $A_i$ is isomorphic to the algebra of $d_i \times d_i$ matrices with
entries in a division algebra $D_i$ over $F$. The action of $A$ on the left
regular representation ${}^\ast\!{A}$ is $a \cdot b = ab$ for $a \in A$, $b \in
{}^\ast\!{A}$, and \eqref{semisimplealgebra} gives the decomposition of
${}^\ast\!{A}$ into isotypic components:
  \begin{equation} \label{semisimplemodule}
  {}^\ast\!{A} = {}^\ast\!{A}_1 \oplus \cdots \oplus {}^\ast\!{A}_r.
  \end{equation}
Each ${}^\ast\!{A}_i$ is the direct sum of $d_i$ isomorphic simple submodules,
and each of these is a $d_i$-dimensional minimal left ideal in $A$ (a column in
the matrix algebra). We consider the direct sum of $t$ copies of ${}^\ast\!{A}$
with the diagonal action:
  \begin{equation} \label{directsum}
  ({}^\ast\!{A})^t
  =
  ( {}^\ast\!{A} )^{[1]} \oplus \cdots \oplus ( {}^\ast\!{A} )^{[t]},
  \qquad
  a \cdot ( b_1, \hdots, b_t ) = ( ab_1, \hdots, ab_t ).
  \end{equation}
If $U \subseteq ({}^\ast\!{A})^t$ is a submodule then usually $U$ is not
homogeneous with respect to \eqref{directsum}:
  \begin{equation} \label{nonhomogeneous}
  U \ne \sum_{k=1}^t{\!\oplus} \Big( U \cap ({}^\ast\!{A})^{[k]} \Big).
  \end{equation}
We combine \eqref{semisimplemodule} and \eqref{directsum} to obtain a finer
decomposition of $({}^\ast\!{A})^t$:
  \begin{equation} \label{isotypic}
  ({}^\ast\!{A})^t
  =
  \sum_{k=1}^t{\!\oplus}
  \big({}^\ast\!{A}\big)^{[k]}
  =
  \sum_{k=1}^t{\!\oplus}
  \sum_{i=1}^r{\!\oplus}
  \big({}^\ast\!{A}\big)^{[k]}_i
  =
  \sum_{i=1}^r{\!\oplus}
  \sum_{k=1}^t{\!\oplus}
  \big({}^\ast\!{A}\big)^{[k]}_i.
  \end{equation}
This gives a direct sum decomposition of $({}^\ast\!{A})^t$ into components
$R_i$:
  \begin{equation} \label{representation}
  ({}^\ast\!{A})^t
  =
  \sum_{i=1}^r{\!\oplus}
  R_i,
  \qquad
  R_i = \sum_{k=1}^t{\!\oplus} \big({}^\ast\!{A}\big)^{[k]}_i.
  \end{equation}
Every submodule $U \subseteq ({}^\ast\!{A})^t$ is homogeneous with respect to
\eqref{representation}.

\begin{lemma}  \label{homogeneouscomponents}
If $U \subseteq ({}^\ast\!{A})^t$ is a submodule then
  \[
  U = \sum_{i=1}^r{\!\oplus} \Big( U \cap R_i \Big).
  \]
\end{lemma}

\begin{proof}
For any $u \in U$ we have $u = u_1 + \cdots + u_r$ where $u_i \in R_i$. If $I_i
\in A_i$ corresponds to the identity matrix then \eqref{semisimplealgebra},
\eqref{directsum}, \eqref{representation} imply $I_i \cdot u = u_i$; hence $u_i
\in U$.
\end{proof}

\subsection{Irreducible representations of the symmetric group}

We apply the general construction to the group algebra $F S_n$ over $F$ of the
symmetric group $S_n$.  We assume that $F = \mathbb{Q}$ or $F = \mathbb{F}_p$
for $p > n$; then $F S_n$ is semisimple by Maschke's theorem. We recall the
structure theory from James and Kerber \cite{JamesKerber}. The irreducible
representations of $S_n$ are in bijection with the partitions of $n$.   Let
$\lambda = (n_1,\hdots,n_\ell)$ be a partition: $n = n_1 + \cdots + n_\ell$
with $n_1 \ge \cdots \ge n_\ell \ge 1$. The frame $[\lambda]$ consists of $n$
empty boxes in $\ell$ left-justified rows with $n_i$ boxes in row $i$.  A
tableau for $\lambda$ is a bijection from $1, \hdots, n$ to the boxes of
$[\lambda]$. In a standard tableau the numbers increase in each row from left
to right and in each column from top to bottom. The number $d_\lambda$ of
standard tableaux with frame $[\lambda]$ is the dimension of the corresponding
irreducible representation.  We have the following direct sum decomposition of
$F S_n$ into orthogonal two-sided ideals isomorphic to simple matrix algebras
over $F$:
  \begin{equation} \label{youngtheorem}
  F S_n \approx \sum_{\lambda}{\!\oplus} A_\lambda,
  \qquad
  A_\lambda = M_{d_\lambda}(F).
  \end{equation}
For us the most important problem is: Given a permutation $\pi$ and a partition
$\lambda$, compute the $d_\lambda \times d_\lambda$ matrix representing $\pi$;
that is, compute the projection of $\pi$ onto the summand $A_\lambda$ in
\eqref{youngtheorem}.  A simple algorithm was found by Clifton \cite{Clifton}
(see also Bergdolt \cite{Bergdolt}). Let $T_1, \hdots, T_d$ ($d = d_\lambda$)
be the standard tableaux for $\lambda$ ordered in some fixed way. Let
$R^\lambda_\pi$ be the matrix defined as follows: \emph{Apply $\pi$ to the
tableau $T_j$. If there exist two numbers that appear together in a column of
$T_i$ and a row of $\pi T_j$, then $(R^\lambda_\pi)_{ij} = 0$.  If not, then
$(R^\lambda_\pi)_{ij}$ equals the sign of the vertical permutation for $T_i$
which leaves the columns of $T_i$ fixed as sets and takes the numbers of $T_i$
into the correct rows they occupy in $\pi T_j$.} The matrix
$R^\lambda_\mathrm{id}$ corresponding to the identity permutation is not
necessarily the identity matrix, but it is always invertible.

\begin{lemma} \cite{Clifton}  \label{cliftonlemma}
The matrix representing $\pi$ in partition $\lambda$ equals
$(R^\lambda_\mathrm{id})^{-1} R^\lambda_\pi$.
\end{lemma}

Since Clifton's algorithm is very important for us, we present it formally in
Figure \ref{cliftonmatrixalgorithm}, following an idea of Hentzel: the
algorithm tries to compute the vertical permutation whose sign gives
$(R^\lambda_\pi)_{ij}$, and returns 0 if it fails.

  \begin{figure}[t]
  \begin{itemize}
  \item Input: A permutation $\pi \in S_n$ and a partition $\lambda =
  (n_1,\dots,n_\ell)$ of $n$.
  \item Output: The Clifton matrix $R^\lambda_\pi$.
  \end{itemize}
  \begin{enumerate}
  \item[(1)]
  Compute the standard tableaux $T_1, \hdots, T_d$ for $\lambda$ where $d = d_\lambda$.
  \item[(2)]
  For $j$ from 1 to $d$ do:
  \begin{enumerate}
  \item[(a)]
  Compute $\pi T_j$.
  \item[(b)]
  For $i$ from 1 to $d$ do:
    \begin{enumerate}
    \item
    Set
    \texttt{ijentry} $\leftarrow$ 1,
    \texttt{number} $\leftarrow$ 1,
    \texttt{finished} $\leftarrow$ \texttt{false}.
    \item
    While $\texttt{number} \le n$ and not \texttt{finished} do:
      \begin{itemize}
      \item
      Set
      \texttt{irow}, \texttt{icol} $\leftarrow$
      row, column indices of \texttt{number} in $T_i$.
      \item
      Set
      \texttt{jrow}, \texttt{jcol} $\leftarrow$
      row, column indices of \texttt{number} in $\pi T_j$.
      \item
      If \texttt{irow} $\ne$ \texttt{jrow} then
        [\texttt{number} \emph{is not in the correct row}]
        \begin{itemize}
        \item
        If \texttt{icol} $> n_{\texttt{jrow}}$ then
        \item[] \quad
        [\emph{the required position does not exist}]
          \item[] \quad set
          \texttt{ijentry} $\leftarrow$ 0,
          \texttt{finished} $\leftarrow$ \texttt{true}
        \item[]
        else if $(T_i)_{\texttt{jrow,icol}} < (T_i)_{\texttt{irow,icol}}$ then
          \item[] \quad
          [\emph{the required position is already occupied}]
          \item[] \quad
          set
          \texttt{ijentry} $\leftarrow$ 0,
          \texttt{finished} $\leftarrow$ \texttt{true}
        \item[]
        else
          \item[] \quad
          [\emph{transpose} \texttt{number} \emph{into the required position}]
          \item[] \quad
          set
          \texttt{ijentry} $\leftarrow$ $-$\texttt{ijentry},
          \item[] \quad
          interchange $(T_i)_{\texttt{irow,icol}}$ and $(T_i)_{\texttt{jrow,icol}}$
        \end{itemize}
      \item
      Set $\texttt{number} \leftarrow \texttt{number} + 1$
      \end{itemize}
      \item Set $(R^\lambda_\pi)_{ij} \leftarrow \texttt{ijentry}$
    \end{enumerate}
  \end{enumerate}
  \item[(3)] Return $R^\lambda_\pi$.
  \end{enumerate}
  \caption{Hentzel's algorithm to compute the Clifton matrix $R^\lambda_\pi$}
  \label{cliftonmatrixalgorithm}
  \end{figure}

\subsection{Polynomial identities and representation theory}

The application of the representation theory of the symmetric group to
polynomial identities was initiated independently by Malcev \cite{Malcev} and
Specht \cite{Specht} in 1950.  The implementation of this theory in computer
algebra was initiated by Hentzel \cite{Hentzel1, Hentzel2} in the 1970's. We
first recall that any polynomial identity (not necessarily multilinear or even
homogeneous) of degree $\le n$ over a field $F$ of characteristic 0 or $p
> n$ is equivalent to a finite set of multilinear identities; see Zhevlakov et
al.~\cite[Chapter 1]{Zhevlakov}. We consider a multilinear nonassociative
identity $I( x_1, \hdots, x_n )$ of degree $n$. We collect the terms which have
the same association type: $I = I_1 + \cdots + I_t$. In each summand $I_k$
every monomial has association type $k$: the monomials differ only by the
permutation of $x_1, \hdots, x_n$. We can therefore regard each $I_k$ as an
element of the group algebra $F S_n$, and the identity $I$ as an element of the
direct sum of $t$ copies of $F S_n$. Following the previous two subsections,
let $U$ be the submodule of $( F S_n )^t$ generated by $I$. Every element of
$U$ is a linear combination of permutations of $I$, and hence is an identity
implied by $I$.  By Lemma \ref{homogeneouscomponents} we know that $U$ is the
direct sum of its components corresponding to the irreducible representations
of $S_n$.  This allows us to study $I$ and its consequences one representation
at a time, so we can break down a large computational problem into much smaller
pieces.

\begin{example}
Consider the Jordan identity $(a^2b)a - a^2(ba)$ in a free commutative
nonassociative algebra.  The multilinear form of this identity (divided by 2)
is
  \[
  u = ((ac)b)d + ((ad)b)c + ((cd)b)a - (ac)(bd) - (ad)(bc) - (cd)(ba).
  \]
In degree 4 there are two association types $((ab)c)d$ and $(ab)(cd)$ for a
commutative nonassociative operation.  We regard $u = u_1 + u_2$ as an element
of the direct sum of two copies of the group algebra $\mathbb{Q} S_4$
corresponding to the two association types:
  \[
  u_1 = ((ac)b)d + ((ad)b)c + ((cd)b)a,
  \qquad
  u_2 = - (ac)(bd) - (ad)(bc) - (cd)(ba).
  \]
To illustrate the inequality of equation \eqref{nonhomogeneous} we note that
the two components $u_1$ and $u_2$ are identities which are not consequences of
the Jordan identity. To illustrate the equality of Lemma
\ref{homogeneouscomponents} we decompose the submodule $U \subseteq (
\mathbb{Q} S_4 )^2$ generated by $u$ into components corresponding to the
irreducible representations of $S_4$. We obtain the linearization of
fourth-power associativity (for $\lambda = 4$) and the identity which says that
the commutator of multiplications is a derivation (for $\lambda = 31$).
\end{example}

\subsection{Ranks and multiplicities}

Let $u^{[1]}, \dots, u^{[g]}$ be a set of multilinear polynomial identities of
degree $n$ over a field $F$ of characteristic 0 or $p > n$.  Suppose that $t$
association types occur in the terms of the $u^{[i]}$, and let $U \subseteq ( F
S_n )^t$ be the submodule generated by $u^{[1]}, \dots, u^{[g]}$. We fix a
partition $\lambda$ of $n$ and write $d = d_\lambda$ for the dimension of the
corresponding irreducible representation of $S_n$.  To determine the
corresponding component of $U$ we construct a $dg \times dt$ matrix $M_\lambda$
with $g$ rows and $t$ columns of $d \times d$ blocks. In block $(i,j)$ we put
the representation matrix for the terms of $u^{[i]}$ in association type $j$,
which can be computed by repeated application of Lemma \ref{cliftonlemma}. We
compute $\mathrm{RCF}(M_\lambda)$, the row canonical form of $M_\lambda$.

\begin{definition} \label{oldrank}
The number of nonzero rows in $\mathrm{RCF}(M_\lambda)$ is the \textbf{rank of
the submodule $U$ in partition $\lambda$}.
\end{definition}

\begin{lemma}
The number of nonzero rows in $\mathrm{RCF}(M_\lambda)$ is the multiplicity of
the irreducible representation corresponding to $\lambda$ in the submodule $U$.
\end{lemma}

  \begin{table}
  \begin{align*}
  &\left[
  \!\!\!
  \begin{array}{ccccc|ccccc}
  \rho_\lambda(E^1_1) & \rho_\lambda(E^1_2) & \cdots &
  \rho_\lambda(E^1_{n-1}) & \rho_\lambda(E^1_n) &
  -I_d   &  O     & \cdots &  O     &  O     \\
  \rho_\lambda(E^2_1) & \rho_\lambda(E^2_2) & \cdots &
  \rho_\lambda(E^2_{n-1}) & \rho_\lambda(E^2_n) &
  O      & -I_d   & \cdots &  O     &  O     \\
  \vdots & \vdots & \ddots & \vdots & \vdots &
  \vdots & \vdots & \ddots & \vdots & \vdots \\
  \rho_\lambda(E^{t-1}_1) & \rho_\lambda(E^{t-1}_2) & \cdots &
  \rho_\lambda(E^{t-1}_{n-1}) & \rho_\lambda(E^{t-1}_n) &
  O      & O      & \cdots & -I_d   &  O     \\
  \rho_\lambda(E^t_1) & \rho_\lambda(E^t_2) & \cdots &
  \rho_\lambda(E^t_{n-1}) & \rho_\lambda(E^t_n) &
  O      & O      & \cdots &  O     & -I_d
  \end{array}
  \!\!\!
  \right]
  \end{align*}
  \caption{Representation matrix in partition $\lambda$ for the dialgebra expansions
  of the right-commutative association types in degree $n$}
  \label{repmatexp}
  \end{table}

We use a modification of this procedure to determine the structure of the
kernel of the expansion map for the quasi-Jordan product.  In degree $n$ there
are $t = R_n$ right-commutative association types (Lemma \ref{typeslemma}) and
$n$ dialgebra association types (corresponding to the position of the center).
We fix a partition $\lambda$ and write $d = d_\lambda$. We create a $td \times
(n{+}t)d$ matrix $X_\lambda$ with $t$ rows and $n{+}t$ columns of $d \times d$
blocks; see Table \ref{repmatexp}. In the right side of $X_\lambda$, in block
$(i,n{+}i)$ for $1 \le i \le t$, we put $-I_d$, the negative of the identity
matrix; the other blocks of the right side are zero. In the left side of
$X_\lambda$, in block $(i,j)$ for $1 \le i \le t$ and $1 \le j \le n$, we put
$\rho_\lambda(E^i_j)$, the representation matrix of the terms with dialgebra
association type $j$ in the expansion of the basic monomial with
right-commutative association type $i$.  The $i$-th row of blocks states that
the basic monomial for the $i$-th right-commutative association type equals its
expansion in the free associative dialgebra. Since the right side of
$X_\lambda$ is the negative of the identity matrix, $X_\lambda$ has full row
rank. We compute $\mathrm{RCF}(X_\lambda)$; there are no zero rows. We
distinguish upper and lower parts of $\mathrm{RCF}(X_\lambda)$: the upper part
contains the rows with leading ones in the left side, and the lower part
contains the rows with leading ones in the right side. The lower left part is
zero; the lower right part represents polynomial identities which are satisfied
by the right-commutative association types as a result of dependence relations
among the dialgebra expansions of the basic right-commutative monomials.

\begin{definition}  \label{allrank}
The number of (nonzero) rows in the lower right block of
$\mathrm{RCF}(X_\lambda)$ will be called the \textbf{rank of all identities
satisfied by the quasi-Jordan product in partition $\lambda$}.
\end{definition}

\begin{lemma}
The number of rows in the lower right block of $\mathrm{RCF}(X_\lambda)$ is the
multiplicity of the irreducible representation corresponding to $\lambda$ in
the kernel of $E_n$.
\end{lemma}

\begin{lemma} \label{newidentitieslemma}
Let the submodule $U$ of lifted identities in degree $n$ for the quasi-Jordan
product be generated by $u^{[i]}, \dots, u^{[g]}$. Let $\lambda$ be a partition
of $n$, let $\mathrm{oldrank}(\lambda)$ be the rank of $U$ in partition
$\lambda$ (Definition \ref{oldrank}), and let $\mathrm{allrank}(\lambda)$ be
the rank of all identities in partition $\lambda$ (Definition \ref{allrank}).
Then $\mathrm{oldrank}(\lambda) \le \mathrm{allrank}(\lambda)$ with equality if
and only if there are no new identities corresponding to partition $\lambda$.
\end{lemma}

\subsection{Rational arithmetic and modular arithmetic}

We prefer to use rational arithmetic, but this is impractical when the matrices
are large: during the computation of the RCF, the numerators and denominators
of the entries can become extremely large, even if the original matrix has
small integer entries. To control the amount of memory required, we use modular
arithmetic with a prime $p$ greater than the degree $n$ of the identities; this
guarantees that $\mathbb{F}_p S_n$ is semisimple. The structure theory of
$\mathbb{Q} S_n$, in particular isomorphism \eqref{youngtheorem}, shows that
the elements of $\mathbb{Q} S_n$ which represent the matrix units in the simple
ideals $A_\lambda$ have coefficients in which the denominators are divisors of
$n!$.  It follows that the $S_n$-module $( F S_n )^t$ has the ``same''
structure over $\mathbb{Q}$ and over $\mathbb{F}_p$ when $p > n$.  Therefore
the ranks we obtain using modular arithmetic will be the same as the ranks we
would have obtained using rational arithmetic. This leaves the problem of
reconstructing rational results from modular results.  In some cases, as in
this paper, modular arithmetic produces coefficients for which the
corresponding rational coefficients are easy to recover: 1, 2, 3, 49, 50, 51,
52, 98, 99, 100 in $\mathbb{F}_{101}$ represent 1, 2, 3, $-3/2$, $-1/2$, $1/2$,
$3/2$, $-3$, $-2$, $-1$ in $\mathbb{Q}$. In other cases, we have to use many
different primes and the Chinese Remainder Theorem; see Bremner and Peresi
\cite{BremnerPeresi}.


\section{Nonexistence of new identities in degree 6: second computation}

  \begin{table}
  \begin{center}
  \begin{tabular}{rlrrrrrrrr}
  & & & \multicolumn{3}{c}{old identities} & \multicolumn{3}{c}{all identities} \\
     & $\lambda$ & $d$ & rows & cols & rank & rows & cols & rank & new \\
   1 & 6         &  1  &   21 &   20 &   17 &   20 &   26 &   17 &   0 \\
   2 & 51        &  5  &  105 &  100 &   85 &  100 &  130 &   85 &   0 \\
   3 & 42        &  9  &  189 &  180 &  153 &  180 &  234 &  153 &   0 \\
   4 & 411       & 10  &  210 &  200 &  172 &  200 &  260 &  172 &   0 \\
   5 & 33        &  5  &  105 &  100 &   85 &  100 &  130 &   85 &   0 \\
   6 & 321       & 16  &  336 &  320 &  274 &  320 &  416 &  274 &   0 \\
   7 & 3111      & 10  &  210 &  200 &  176 &  200 &  260 &  176 &   0 \\
   8 & 222       &  5  &  105 &  100 &   85 &  100 &  130 &   85 &   0 \\
   9 & 2211      &  9  &  189 &  180 &  157 &  180 &  234 &  157 &   0 \\
  10 & 21111     &  5  &  105 &  100 &   91 &  100 &  130 &   91 &   0 \\
  11 & 111111    &  1  &   21 &   20 &   19 &   20 &   26 &   19 &   0
  \end{tabular}
  \end{center}
  \caption{Degree 6: matrix ranks for all representations}
  \label{degree6ranks}
  \end{table}

Table \ref{degree6ranks} gives the matrix ranks for $p = 101$. There are no new
identities, confirming our earlier computations without the representation
theory of $S_n$.

When we use representation theory, we have two kinds of lifted identities: the
31 symmetries of the right-commutative association types and the 56 liftings of
the identities $J$ and $K$. In degree 6, we have 20 association types. For
partition $\lambda$ with dimension $d = d_\lambda$, we create a matrix with
$20d$ columns and $21d$ rows, initialized to zero. For each lifted identity we
put the representation matrices in the bottom $d$ rows and compute the RCF. In
Table \ref{degree6ranks} under ``old identities'' the columns ``rows'' and
``cols'' contain $21d$ and $20d$; the column ``rank'' contains the final rank
of the matrix. From the complete list of 56 identities obtained by lifting $J$
and $K$ to degree 6, we retain only those identities which increase the rank
for at least one partition. We recover the same 25 generators that we obtained
earlier.

To compute all the identities for partition $\lambda$, we create a matrix with
$20d$ rows and $26d$ columns, with a $20d \times 6d$ left block corresponding
to the dialgebra expansions, and a $20d \times 20d$ right block corresponding
to the basic right-commutative monomials. In Table \ref{degree6ranks} under
``all identities'' the columns ``rows'' and ``cols'' contain $20d$ and $26d$;
the column ``rank'' contains the number of rows in the RCF which have leading
ones in the lower right block of the matrix: these rows represent polynomial
identities satisfied by the quasi-Jordan product.

When the two ranks are the same for partition $\lambda$, there are no new
identities for this representation (Lemma \ref{newidentitieslemma}). We checked
this by verifying that the two matrices are in fact equal. Let $r$ be the
common rank for partition $\lambda$. The first matrix has size $r \times 20d$;
these are the nonzero rows of the RCF of the matrix for the lifted identities.
The second matrix has the same size; it contains the rows of the RCF of the
matrix for the expansion identities with leading ones in the lower right block.


\section{Nonexistence of new identities in degree 7}

  \begin{table}
  \begin{center}
  \begin{tabular}{rlrrrrrrrr}
  & &
  & \multicolumn{3}{c}{old identities}
  & \multicolumn{3}{c}{all identities}
  \\
     & $\lambda$  & $d$ & rows & cols & rank & rows & cols & rank & new \\
   1 & 7          &  1  &   47 &   46 &   42 &   46 &   53 &   42 &   0 \\
   2 & 61         &  6  &  282 &  276 &  255 &  276 &  318 &  255 &   0 \\
   3 & 52         & 14  &  658 &  644 &  594 &  644 &  742 &  594 &   0 \\
   4 & 511        & 15  &  705 &  690 &  641 &  690 &  795 &  641 &   0 \\
   5 & 43         & 14  &  658 &  644 &  595 &  644 &  742 &  595 &   0 \\
   6 & 421        & 35  & 1645 & 1610 & 1490 & 1610 & 1855 & 1490 &   0 \\
   7 & 4111       & 20  &  940 &  920 &  859 &  920 & 1060 &  859 &   0 \\
   8 & 331        & 21  &  987 &  966 &  895 &  966 & 1113 &  895 &   0 \\
   9 & 322        & 21  &  987 &  966 &  892 &  966 & 1113 &  892 &   0 \\
  10 & 3211       & 35  & 1645 & 1610 & 1499 & 1610 & 1855 & 1499 &   0 \\
  11 & 31111      & 15  &  705 &  690 &  651 &  690 &  795 &  651 &   0 \\
  12 & 2221       & 14  &  658 &  644 &  598 &  644 &  742 &  598 &   0 \\
  13 & 22111      & 14  &  658 &  644 &  607 &  644 &  742 &  607 &   0 \\
  14 & 211111     &  6  &  282 &  276 &  265 &  276 &  318 &  265 &   0 \\
  15 & 1111111    &  1  &   47 &   46 &   45 &   46 &   53 &   45 &   0
  \end{tabular}
  \end{center}
  \caption{Degree 7: matrix ranks for all representations}
  \label{degree7ranks}
  \end{table}

Table \ref{degree7ranks} gives the ranks for $p = 101$: there are no new
identities in degree 7.  The computations are similar to those for degree 6,
except that the matrices are larger.  The lifted (``old'') identities consist
of 89 symmetries of the 46 right-commutative association types and 200 liftings
of $J$ and $K$. The matrix of ``old identities'' has size $47d \times 46d$, and
the matrix of ``all identities'' has size $46d \times 53d$. A subset of 55
identities suffices to generate all the lifted identities.


\section{Special identities} \label{sectionspecial}

The previous sections show that there are no special identities of degree $\le
7$ for quasi-Jordan algebras.  The referee suggested the basic ideas for the
following algorithm to construct certain special identities.

\subsection{Construction of noncommutative preimages of the Glennie identity}
\label{referee'salgorithm}

Recall the Jordan triple product $\{abc\} = (ab)c + (cb)a - (ac)b$ and the
Glennie identity of degree 8 from \cite[page 79]{Zhevlakov}:
  \begin{align*}
  G
  &=
  2\{\{b\{aca\}b\}c(ab)\}
  -
  \{b\{a\{c(ab)c\}a\}b\}
  \\
  &\quad
  -
  2\{(ab)c\{a\{bcb\}a\}\}
  +
  \{a\{b\{c(ab)c\}b\}a\}.
  \end{align*}
We have the following algorithm for constructing noncommutative preimages of
the Glennie identity:
  \begin{enumerate}
  \item Start with the Glennie identity $G$, a commutative nonassociative
      polynomial, homogeneous of degree 8, with variables $aaabbbcc$.
      Choose either $x = a$ or $x = b$ or $x = c$ as the variable to
      linearize.
  \item Partially linearize $G$ by replacing each occurrence of $x$ in each
      term by the new variable $d$ in every possible position; that is,
      apply the operator $\Delta^1_x(d)$ from \cite[Chapter 1]{Zhevlakov}.
      For example, if $x = a$ and a term has the form ${-}a{-}a{-}a{-}$
      then the term becomes the sum of the three new terms
      \[
      {-}d{-}a{-}a{-}, \quad {-}a{-}d{-}a{-}, \quad {-}a{-}a{-}d{-}.
      \]
      We obtain a special Jordan identity $H$ which is linear in $d$.
  \item Apply the algorithm of Pozhidaev \cite[Section 3]{Pozhidaev}:
      convert each algebra monomial into a dialgebra monomial by making $d$
      the center. Operations to the left of $d$ become $\vdash$ and
      operations to the right of $d$ become $\dashv$. We obtain a
      nonassociative dialgebra polynomial $I$.
  \item Replace the right dialgebra operation $\vdash$ in $I$ by the
      opposite of the left dialgebra operation $\dashv$; that is, replace
      every occurrence of $y \vdash z$ by $z \dashv y$. (Recall that $y
      \vdash z = z \dashv y$ is the dialgebra consequence of
      commutativity.) Remove the symbols $\dashv$ from $I$. We obtain a
      noncommutative nonassociative polynomial $J$ in one binary operation
      written as juxtaposition.
  \item Replace every occurrence of $d$ in $J$ by the original variable
      $x$. We obtain a nonassociative polynomial $K$, homogeneous of degree
      8, with variables $aaabbbcc$; this is a noncommutative preimage of
      the Glennie identity $G$.
  \item Expand $K$ into the free associative dialgebra using the
      quasi-Jordan product: that is, each product $yz$ in $K$ becomes $y
      \dashv z + z \vdash y$. We convert each term in the expansion to its
      normal form, collect similar terms, and verify that the result
      collapses to zero in all cases: $x = a$, $x = b$, and $x = c$.
  \end{enumerate}
For $x = a, b$ the identity $K$ has 100 terms. For $x = c$ the identity $K$ has
72 terms; see Table \ref{preimage}. (The monomials have been sorted and the
coefficients divided by $-2$.) It is not clear \emph{a priori} that the
preimages obtained by this algorithm hold in any special quasi-Jordan algebra,
but we have verified this by direct computation using Maple in step (6) of the
algorithm.

\subsection{A new special identity}

The preceding arguments prove that special identities exist and that some of
them are noncommutative preimages of the Glennie identity. We have discovered a
new special identity, homogeneous of degree 8, with variables $aaaabbbc$; this
identity has three variables and is linear in one variable. If it was a
noncommutative preimage of the Glennie identity, then its commutative version
would be satisfied by all special Jordan algebras and hence by all Jordan
algebras (Macdonald's theorem); but its commutative version would be the
Glennie identity, and this is a contradiction.

  \begin{table}
  \[
  \begin{array}{rrr}
   4 (((((c(ab))c)a)a)b)b &  
  -4 (((((c(ab))c)b)b)a)a &  
  -4 (((((cc)(ab))a)a)b)b \\ 
  +4 (((((cc)(ab))b)b)a)a &  
  +4 ((((c((ba)c))a)a)b)b &  
  -4 ((((c((ba)c))b)b)a)a \\ 
  -2 ((((c(ab))c)(aa))b)b &  
  +2 ((((c(ab))c)(bb))a)a &  
  +2 ((((cc)(ab))(aa))b)b \\ 
  -2 ((((cc)(ab))(bb))a)a &  
  -2 (((c((ba)c))(aa))b)b &  
  +2 (((c((ba)c))(bb))a)a \\ 
  +4 (((((ca)a)b)b)(ab))c &  
  -4 (((((cb)b)a)a)(ba))c &  
  -2 ((((c(aa))b)b)(ab))c \\ 
  +2 ((((c(bb))a)a)(ba))c &  
  -2 ((((ca)a)(bb))(ab))c &  
  +2 ((((cb)b)(aa))(ba))c \\ 
  +  (((c(aa))(bb))(ab))c &  
  -  (((c(bb))(aa))(ba))c &  
  -4 (((((ca)a)b)b)c)(ab) \\ 
  +4 (((((cb)b)a)a)c)(ab) &  
  +2 ((((c(aa))b)b)c)(ab) &  
  -2 ((((c(ab))c)a)a)(bb) \\ 
  +2 ((((c(ab))c)b)b)(aa) &  
  -2 ((((c(bb))a)a)c)(ab) &  
  +2 ((((cc)(ab))a)a)(bb) \\ 
  -2 ((((cc)(ab))b)b)(aa) &  
  -2 (((c((ba)c))a)a)(bb) &  
  +2 (((c((ba)c))b)b)(aa) \\ 
  +2 ((((ca)a)(bb))c)(ab) &  
  -2 ((((cb)b)(aa))c)(ab) &  
  -  (((c(aa))(bb))c)(ab) \\ 
  +  (((c(bb))(aa))c)(ab) &  
  +  (((c(ab))c)(aa))(bb) &  
  -  (((c(ab))c)(bb))(aa) \\ 
  -  (((cc)(ab))(aa))(bb) &  
  +  (((cc)(ab))(bb))(aa) &  
  +  ((c((ba)c))(aa))(bb) \\ 
  -  ((c((ba)c))(bb))(aa) &  
  +2 (c((a(ca))(bb)))(ab) &  
  -2 (c((b(cb))(aa)))(ab) \\ 
  -  (c((c(aa))(bb)))(ab) &  
  +  (c((c(bb))(aa)))(ab) &  
  +4 (c(a((b(cb))a)))(ab) \\ 
  -2 (c(a((c(bb))a)))(ab) &  
  -4 (c(b((a(ca))b)))(ab) &  
  +2 (c(b((c(aa))b)))(ab) \\ 
  -4 ((((ca)a)b)b)(c(ba)) &  
  +4 ((((cb)b)a)a)(c(ba)) &  
  +2 (((c(aa))b)b)(c(ba)) \\ 
  -2 (((c(bb))a)a)(c(ba)) &  
  +2 (((ca)a)(bb))(c(ba)) &  
  -2 (((cb)b)(aa))(c(ba)) \\ 
  -  ((c(aa))(bb))(c(ba)) &  
  +  ((c(bb))(aa))(c(ba)) &  
  -2 (c(ba))((a((bb)c))a) \\ 
  +4 (c(ba))((a((bc)b))a) &  
  +2 (c(ba))((b((aa)c))b) &  
  -4 (c(ba))((b((ac)a))b) \\ 
  +  (c(ba))((aa)((bb)c)) &  
  -2 (c(ba))((aa)((bc)b)) &  
  -  (c(ba))((bb)((aa)c)) \\ 
  +2 (c(ba))((bb)((ac)a)) &  
  +2 c(((b(cb))(aa))(ba)) &  
  -  c(((c(bb))(aa))(ba)) \\ 
  -4 c((a((b(cb))a))(ba)) &  
  +2 c((a((c(bb))a))(ba)) &  
  -2 c((ba)((a(ca))(bb))) \\ 
  +  c((ba)((c(aa))(bb))) &  
  +4 c((ba)(b((a(ca))b))) &  
  -2 c((ba)(b((c(aa))b))) \\ 
  \end{array}
  \]
  \caption{Noncommutative preimage of the Glennie identity ($x = c$)}
  \label{preimage}
  \end{table}


\section{New identities in degree 8}

  \begin{table}
  \begin{center}
  \begin{tabular}{rlrrrrrrrrr}
  & & & \multicolumn{3}{c}{old identities} & \multicolumn{3}{c}{all identities} \\
     & $\lambda$ & $d$ & rows & cols & rank & rows &  cols & rank & new \\
   1 & 8         &  1  &  107 &  106 &  102 &  106 &   114 &  102 &   0 \\
   2 & 71        &  7  &  749 &  742 &  714 &  742 &   798 &  714 &   0 \\
   3 & 62        & 20  & 2140 & 2120 & 2040 & 2120 &  2280 & 2040 &   0 \\
   4 & 611       & 21  & 2247 & 2226 & 2145 & 2226 &  2394 & 2145 &   0 \\
   5 & 53        & 28  & 2996 & 2968 & 2856 & 2968 &  3192 & 2856 &   0 \\
   6 & 521       & 64  & 6848 & 6784 & 6532 & 6784 &  7296 & 6532 &   0 \\
   7 & 5111      & 35  & 3745 & 3710 & 3582 & 3710 &  3990 & 3582 &   0 \\
   8 & 44        & 14  & 1498 & 1484 & 1428 & 1484 &  1596 & 1428 &   0 \\
   9 & 431       & 70  & 7490 & 7420 & 7142 & 7420 &  7980 & 7143 &   1 & $\leftarrow$ \\
  10 & 422       & 56  & 5992 & 5936 & 5712 & 5936 &  6384 & 5713 &   1 & $\leftarrow$ \\
  11 & 4211      & 90  & 9630 & 9540 & 9199 & 9540 & 10260 & 9199 &   0 \\
  12 & 41111     & 35  & 3745 & 3710 & 3594 & 3710 &  3990 & 3594 &   0 \\
  13 & 332       & 42  & 4494 & 4452 & 4284 & 4452 &  4788 & 4286 &   2 & $\leftarrow$ \\
  14 & 3311      & 56  & 5992 & 5936 & 5722 & 5936 &  6384 & 5723 &   1 & $\leftarrow$ \\
  15 & 3221      & 70  & 7490 & 7420 & 7149 & 7420 &  7980 & 7150 &   1 & $\leftarrow$ \\
  16 & 32111     & 64  & 6848 & 6784 & 6565 & 6784 &  7296 & 6565 &   0 \\
  17 & 311111    & 21  & 2247 & 2226 & 2169 & 2226 &  2394 & 2169 &   0 \\
  18 & 2222      & 14  & 1498 & 1484 & 1429 & 1484 &  1596 & 1429 &   0 \\
  19 & 22211     & 28  & 2996 & 2968 & 2870 & 2968 &  3192 & 2870 &   0 \\
  20 & 221111    & 20  & 2140 & 2120 & 2065 & 2120 &  2280 & 2065 &   0 \\
  21 & 2111111   &  7  &  749 &  742 &  729 &  742 &   798 &  729 &   0 \\
  22 & 11111111  &  1  &  107 &  106 &  105 &  106 &   114 &  105 &   0
  \end{tabular}
  \end{center}
  \caption{Degree 8: matrix ranks for all representations}
  \label{degree8ranks}
  \end{table}

The computations in degree 8 are similar to degree 7, except that the matrices
are larger; see Table \ref{degree8ranks}. There are no new identities except in
the representations corresponding to $\lambda = 431, 422, 332, 3311, 3221$
where the difference between allrank($\lambda$) and oldrank($\lambda$) is 1, 1,
2, 1, 1 respectively. The lifted identities consist of 242 symmetries of the
106 association types together with 495 liftings of $J$ and $K$. The matrix of
``old identities'' has size $107d \times 106d$, and the matrix of ``all
identities'' has size $106d \times 114d$. The ranks were computed using modular
arithmetic with $p = 101$. We can recover rational results from modular
results, and use rational arithmetic to verify the results; see the next
section for the case $\lambda = 431$.

\begin{definition}
We say that a polynomial identity in degree $n$ is \textbf{irreducible} if its
complete linearization generates an irreducible representation of $S_n$.
\end{definition}

\begin{theorem}
There are six new irreducible identities for the quasi-Jordan product in degree
8: one each for $\lambda = 431, 422, 3311, 3221$ and two for $\lambda = 332$.
\end{theorem}

From this we obtain the following example of an exceptional (non-special)
quasi-Jordan algebra.

\begin{corollary}
There exists a noncommutative nilpotent exceptional quasi-Jordan algebra.
\end{corollary}

\begin{proof}
Let $\mathcal{N}$ be the variety of quasi-Jordan algebras defined by the
identities $x_1 \cdots x_9 = 0$ with any placement of parentheses. Let $X$ be
the free algebra in $\mathcal{N}$ on the generators $a, b, c$.  The
quasi-Jordan polynomial discussed in the next section is nonzero in $X$ but is
zero in every special quasi-Jordan algebra.
\end{proof}


\section{A special identity for partition 431} \label{section431}

  \begin{table}
  \begin{center}
  \[
  \begin{array}{rrr}
   2 ((((((aa)a)a)b)b)c)b &
  -2 ((((((aa)a)b)a)b)c)b &
  +2 ((((((aa)a)b)b)a)b)c \\[-1pt]
  -2 ((((((aa)a)b)b)b)a)c &
  +2 ((((((aa)a)b)b)b)c)a &
  -2 ((((((aa)a)b)b)c)a)b \\[-1pt]
  -2 ((((((aa)a)b)b)c)b)a &
  +4 ((((((aa)a)b)c)b)b)a &
  +2 ((((((aa)a)c)b)a)b)b \\[-1pt]
  -4 ((((((aa)a)c)b)b)a)b &
  -2 ((((((aa)b)a)a)b)b)c &
  +4 ((((((aa)b)a)b)a)b)c \\[-1pt]
  -2 ((((((aa)b)a)b)b)a)c &
  -2 ((((((aa)b)a)b)c)a)b &
  +2 ((((((aa)b)a)b)c)b)a \\[-1pt]
  +2 ((((((aa)b)a)c)a)b)b &
  -2 ((((((aa)b)a)c)b)a)b &
  +2 ((((((aa)b)b)c)b)a)a \\[-1pt]
  -2 ((((((aa)b)c)b)a)a)b &
  +2 ((((((aa)b)c)b)a)b)a &
  -2 ((((((aa)b)c)b)b)a)a \\[-1pt]
  +2 ((((((ab)a)a)a)b)b)c &
  -2 ((((((ab)a)a)b)a)b)c &
  -2 ((((((ab)a)a)b)a)c)b \\[-1pt]
  -2 ((((((ab)a)a)b)c)b)a &
  -2 ((((((ab)a)a)c)a)b)b &
  +2 ((((((ab)a)a)c)b)a)b \\[-1pt]
  -4 ((((((ab)a)b)a)a)b)c &
  +2 ((((((ab)a)b)a)a)c)b &
  +4 ((((((ab)a)b)a)b)a)c \\[-1pt]
  +4 ((((((ab)a)b)a)c)a)b &
  -2 ((((((ab)a)b)b)a)c)a &
  +2 ((((((ab)a)b)c)a)a)b \\[-1pt]
  -6 ((((((ab)a)b)c)a)b)a &
  -2 ((((((ab)a)c)a)a)b)b &
  +4 ((((((ab)a)c)a)b)a)b \\[-1pt]
  -2 ((((((ab)a)c)a)b)b)a &
  +2 ((((((ab)a)c)b)a)a)b &
  +2 ((((((ab)a)c)b)b)a)a \\[-1pt]
  +2 ((((((ab)b)a)b)c)a)a &
  -2 ((((((ab)b)a)c)a)a)b &
  +2 ((((((ab)b)a)c)a)b)a \\[-1pt]
  -2 ((((((ab)b)a)c)b)a)a &
  +2 ((((((ab)b)c)a)b)a)a &
  -2 ((((((ab)b)c)b)a)a)a \\[-1pt]
  -2 ((((((ab)c)a)b)a)a)b &
  +2 ((((((ab)c)a)b)a)b)a &
  -2 ((((((ab)c)a)b)b)a)a \\[-1pt]
  +2 ((((((ab)c)b)a)a)a)b &
  -2 ((((((ab)c)b)a)a)b)a &
  +2 ((((((ab)c)b)b)a)a)a \\[-1pt]
  +2 ((((((ac)a)a)b)a)b)b &
  -2 ((((((ac)a)a)b)b)a)b &
  +2 ((((((ac)a)a)b)b)b)a \\[-1pt]
  -2 ((((((ac)a)b)a)a)b)b &
  +2 ((((((ac)a)b)a)b)a)b &
  -4 ((((((ac)a)b)a)b)b)a \\[-1pt]
  +2 ((((((ac)a)b)b)a)a)b &
  +2 ((((((ac)a)b)b)a)b)a &
  -2 ((((((ac)a)b)b)b)a)a \\[-1pt]
  -2 ((((((ac)b)a)b)a)a)b &
  +2 ((((((ac)b)a)b)a)b)a &
  +2 ((((((ac)b)a)b)b)a)a \\[-1pt]
  -2 ((((((ac)b)b)a)b)a)a &
  -2 (((((a(aa))a)b)b)c)b &
  +2 (((((a(aa))b)a)b)c)b \\[-1pt]
  -2 (((((a(aa))b)b)a)b)c &
  +2 (((((a(aa))b)b)b)a)c &
  -2 (((((a(aa))b)b)b)c)a \\[-1pt]
  +2 (((((a(aa))b)b)c)a)b &
  +2 (((((a(aa))b)b)c)b)a &
  -4 (((((a(aa))b)c)b)b)a \\[-1pt]
  -2 (((((a(aa))c)b)a)b)b &
  +4 (((((a(aa))c)b)b)a)b &
  -2 (((((a(ab))a)a)b)b)c \\[-1pt]
  -2 (((((a(ab))a)a)c)b)b &
  +2 (((((a(ab))a)b)a)b)c &
  +4 (((((a(ab))a)b)a)c)b \\[-1pt]
  -2 (((((a(ab))a)b)c)a)b &
  +2 (((((a(ab))a)b)c)b)a &
  +4 (((((a(ab))a)c)a)b)b \\[-1pt]
  -4 (((((a(ab))a)c)b)a)b &
  +2 (((((a(ab))a)c)b)b)a &
  +4 (((((a(ab))b)a)a)b)c \\[-1pt]
  -6 (((((a(ab))b)a)b)a)c &
  +2 (((((a(ab))b)a)b)c)a &
  -4 (((((a(ab))b)a)c)a)b \\[-1pt]
  -2 (((((a(ab))b)a)c)b)a &
  +2 (((((a(ab))b)b)a)c)a &
  -2 (((((a(ab))b)b)c)a)a \\[-1pt]
  -2 (((((a(ab))b)c)a)a)b &
  +6 (((((a(ab))b)c)a)b)a &
  +2 (((((a(ab))c)a)a)b)b \\[-1pt]
  -2 (((((a(ab))c)a)b)a)b &
  -2 (((((a(ab))c)b)a)b)a &
  -4 (((((a(ac))a)b)a)b)b \\[-1pt]
  +4 (((((a(ac))a)b)b)a)b &
  -2 (((((a(ac))a)b)b)b)a &
  +2 (((((a(ac))b)a)a)b)b \\[-1pt]
  -2 (((((a(ac))b)b)a)b)a &
  +2 (((((a(ac))b)b)b)a)a &
  +  (((((a(bb))a)a)a)c)b \\[-1pt]
  -  (((((a(bb))a)a)c)a)b &
  -  (((((a(bb))a)b)a)a)c &
  +  (((((a(bb))a)b)a)c)a \\[-1pt]
  +  (((((a(bb))a)b)c)a)a &
  -  (((((a(bb))a)c)b)a)a &
  -  (((((a(bb))b)c)a)a)a \\[-1pt]
  +  (((((a(bb))c)b)a)a)a &
  +2 (((((a(bc))a)b)a)a)b &
  -2 (((((a(bc))a)b)a)b)a \\[-1pt]
  +2 (((((a(bc))b)a)b)a)a &
  +2 (((((aa)(ab))a)b)b)c &
  -2 (((((aa)(ab))a)b)c)b \\[-1pt]
  +4 (((((aa)(ab))a)c)b)b &
  -2 (((((aa)(ab))b)a)b)c &
  -4 (((((aa)(ab))b)a)c)b \\[-1pt]
  +2 (((((aa)(ab))b)b)a)c &
  -2 (((((aa)(ab))b)b)c)a &
  +6 (((((aa)(ab))b)c)a)b \\[-1pt]
  -4 (((((aa)(ab))c)a)b)b &
  +4 (((((aa)(ab))c)b)a)b &
  -2 (((((aa)(ab))c)b)b)a \\[-1pt]
  -3 (((((aa)(bb))a)a)c)b &
  +2 (((((aa)(bb))a)b)a)c &
  -2 (((((aa)(bb))a)b)c)a \\[-1pt]
  +2 (((((aa)(bb))a)c)a)b &
  +2 (((((aa)(bb))a)c)b)a &
  -  (((((aa)(bb))b)a)c)a \\[-1pt]
  +2 (((((aa)(bb))c)a)a)b &
  -2 (((((aa)(bb))c)a)b)a &
  +  (((((aa)(bb))c)b)a)a \\[-1pt]
  +2 (((((ab)(ac))a)a)b)b &
  -2 (((((ab)(ac))a)b)a)b &
  +2 (((((ab)(ac))a)b)b)a \\[-1pt]
  -2 (((((ab)(ac))b)a)a)b &
  +2 (((((ab)(ac))b)a)b)a &
  -2 (((((ab)(ac))b)b)a)a \\[-1pt]
  -2 (((((ab)(bc))a)a)a)b &
  +2 (((((ab)(bc))a)a)b)a &
  -2 (((((ab)(bc))a)b)a)a \\[-1pt]
  -2 (((((ab)(bc))b)a)a)a &
  +2 (((((aa)a)(ab))b)c)b &
  -2 (((((aa)a)(ab))c)b)b \\[-1pt]
  -3 (((((aa)a)(bb))a)b)c &
  +4 (((((aa)a)(bb))a)c)b &
  +2 (((((aa)a)(bb))b)a)c \\[-1pt]
  +  (((((aa)a)(bb))c)a)b &
  -2 (((((aa)a)(bb))c)b)a &
  -2 (((((aa)a)(bc))a)b)b
  \end{array}
  \]
  \end{center}
  \caption{Special identity for partition 431 (terms 1 to 150)}
  \label{specialidentitypart1}
  \end{table}

  \begin{table}
  \begin{center}
  \[
  \begin{array}{rrr}
  +4 (((((aa)a)(bc))b)a)b &
  +4 (((((aa)b)(ac))a)b)b &
  +2 (((((aa)b)(ac))b)a)b \\[-1pt]
  +2 (((((aa)b)(ac))b)b)a &
  +2 (((((aa)b)(bc))a)a)b &
  +2 (((((aa)b)(bc))a)b)a \\[-1pt]
  +2 (((((ab)a)(ab))a)b)c &
  -2 (((((ab)a)(ab))a)c)b &
  -2 (((((ab)a)(ab))b)a)c \\[-1pt]
  +2 (((((ab)a)(ab))b)c)a &
  -6 (((((ab)a)(ab))c)a)b &
  -4 (((((ab)a)(ac))a)b)b \\[-1pt]
  -6 (((((ab)a)(ac))b)a)b &
  -2 (((((ab)a)(ac))b)b)a &
  +  (((((ab)a)(bb))a)c)a \\[-1pt]
  -2 (((((ab)a)(bb))c)a)a &
  -4 (((((ab)a)(bc))a)a)b &
  -2 (((((ab)a)(bc))a)b)a \\[-1pt]
  -2 (((((ab)a)(bc))b)a)a &
  -2 (((((ab)b)(ac))a)b)a &
  +2 (((((ab)b)(ac))b)a)a \\[-1pt]
  +2 (((((ab)b)(bc))a)a)a &
  -4 (((((ac)a)(ab))b)a)b &
  +2 (((((ac)a)(ab))b)b)a \\[-1pt]
  -  (((((ac)a)(bb))a)a)b &
  -  (((((ac)a)(bb))a)b)a &
  +  (((((ac)a)(bb))b)a)a \\[-1pt]
  +2 ((((a(aa))(bb))a)b)c &
  -  ((((a(aa))(bb))a)c)b &
  -2 ((((a(aa))(bb))b)a)c \\[-1pt]
  +2 ((((a(aa))(bb))b)c)a &
  -2 ((((a(aa))(bb))c)a)b &
  +2 ((((a(aa))(bc))a)b)b \\[-1pt]
  -4 ((((a(aa))(bc))b)a)b &
  -2 ((((a(ab))(ab))a)b)c &
  +4 ((((a(ab))(ab))b)a)c \\[-1pt]
  -4 ((((a(ab))(ab))b)c)a &
  +4 ((((a(ab))(ab))c)a)b &
  +4 ((((a(ab))(ac))b)a)b \\[-1pt]
  +  ((((a(ab))(bb))a)a)c &
  -2 ((((a(ab))(bb))a)c)a &
  +  ((((a(ab))(bb))c)a)a \\[-1pt]
  +2 ((((a(ab))(bc))a)a)b &
  +2 ((((a(ab))(bc))b)a)a &
  +  ((((a(ac))(bb))a)b)a \\[-1pt]
  -  ((((a(ac))(bb))b)a)a &
  +  (((((aa)a)a)(bb))b)c &
  -2 (((((aa)a)a)(bb))c)b \\[-1pt]
  +2 (((((aa)a)b)(ac))b)b &
  -2 (((((aa)a)b)(bc))b)a &
  -2 (((((aa)b)a)(ab))b)c \\[-1pt]
  +2 (((((aa)b)a)(ab))c)b &
  -4 (((((aa)b)a)(ac))b)b &
  -  (((((aa)b)a)(bb))a)c \\[-1pt]
  -4 (((((aa)b)a)(bc))a)b &
  -2 (((((aa)b)a)(bc))b)a &
  -6 (((((aa)b)b)(ac))a)b \\[-1pt]
  -4 (((((aa)b)b)(ac))b)a &
  -2 (((((aa)b)b)(bc))a)a &
  -2 (((((aa)c)a)(bb))a)b \\[-1pt]
  +  (((((aa)c)a)(bb))b)a &
  +2 (((((ab)a)a)(ac))b)b &
  +  (((((ab)a)a)(bb))a)c \\[-1pt]
  -2 (((((ab)a)a)(bb))c)a &
  +4 (((((ab)a)a)(bc))a)b &
  +6 (((((ab)a)b)(ac))a)b \\[-1pt]
  +6 (((((ab)a)b)(ac))b)a &
  +4 (((((ab)a)b)(bc))a)a &
  +2 (((((ab)b)a)(ab))c)a \\[-1pt]
  +4 (((((ab)b)a)(ac))a)b &
  +2 (((((ab)b)a)(ac))b)a &
  +2 (((((ab)b)a)(bc))a)a \\[-1pt]
  +4 (((((ab)c)a)(ab))a)b &
  +2 (((((ab)c)a)(bb))a)a &
  +  (((((ac)a)a)(bb))a)b \\[-1pt]
  -  (((((ac)a)a)(bb))b)a &
  -  (((((ac)b)a)(bb))a)a &
  +  ((((a(aa))a)(bb))c)b \\[-1pt]
  +4 ((((a(aa))b)(bc))b)a &
  -2 ((((a(ab))a)(bb))a)c &
  +3 ((((a(ab))a)(bb))c)a \\[-1pt]
  -4 ((((a(ab))b)(ac))a)b &
  -4 ((((a(ab))b)(ac))b)a &
  -4 ((((a(ab))b)(bc))a)a \\[-1pt]
  +2 (((((aa)a)a)b)(bc))b &
  -2 (((((aa)a)b)a)(bb))c &
  -4 (((((aa)a)b)a)(bc))b \\[-1pt]
  -2 (((((aa)a)b)b)(ac))b &
  +2 (((((aa)b)a)a)(bb))c &
  +2 (((((aa)b)a)a)(bc))b \\[-1pt]
  +2 (((((aa)b)a)b)(bc))a &
  +4 (((((aa)b)b)a)(ac))b &
  +2 (((((aa)b)b)a)(bc))a \\[-1pt]
  +2 (((((aa)b)b)b)(ac))a &
  -2 (((((aa)b)c)a)(ab))b &
  -  (((((ab)a)a)a)(bb))c \\[-1pt]
  +2 (((((ab)a)a)b)(ac))b &
  +4 (((((ab)a)a)b)(bc))a &
  +2 (((((ab)a)b)a)(ab))c \\[-1pt]
  -2 (((((ab)a)b)a)(ac))b &
  -4 (((((ab)a)b)a)(bc))a &
  +4 (((((ab)a)c)a)(ab))b \\[-1pt]
  +2 (((((ab)a)c)a)(bb))a &
  -2 (((((ab)b)a)a)(ac))b &
  -2 (((((ab)b)a)a)(bc))a \\[-1pt]
  -6 (((((ab)b)a)b)(ac))a &
  -2 (((((ab)b)c)a)(ab))a &
  -2 (((((ab)c)a)a)(ab))b \\[-1pt]
  -  (((((ab)c)a)a)(bb))a &
  +2 (((((ac)a)b)a)(ab))b &
  +  (((((ac)a)b)a)(bb))a \\[-1pt]
  -2 ((((a(aa))b)b)(bc))a &
  +  ((((a(ab))a)a)(bb))c &
  -2 ((((a(ab))a)b)(ac))b \\[-1pt]
  -4 ((((a(ab))a)b)(bc))a &
  -2 ((((a(ab))b)a)(ab))c &
  +2 ((((a(ab))b)a)(ac))b \\[-1pt]
  +4 ((((a(ab))b)a)(bc))a &
  +4 ((((a(ab))b)b)(ac))a &
  -2 ((((a(ab))c)a)(ab))b \\[-1pt]
  -  ((((a(ab))c)a)(bb))a &
  -2 (((((aa)a)a)b)b)(bc) &
  +4 (((((aa)a)b)a)b)(bc) \\[-1pt]
  -2 (((((aa)b)a)b)a)(bc) &
  +2 (((((aa)b)a)b)b)(ac) &
  -2 (((((aa)b)b)a)b)(ac) \\[-1pt]
  +2 (((((aa)b)b)c)a)(ab) &
  -  (((((aa)b)c)a)a)(bb) &
  -2 (((((aa)b)c)b)a)(ab) \\[-1pt]
  +  (((((aa)c)a)b)a)(bb) &
  +  (((((ab)a)a)c)a)(bb) &
  -2 (((((ab)a)b)a)b)(ac) \\[-1pt]
  -2 (((((ab)a)b)c)a)(ab) &
  +2 (((((ab)a)c)b)a)(ab) &
  +2 (((((ab)b)a)b)a)(ac) \\[-1pt]
  -2 (((((ab)b)a)c)a)(ab) &
  +2 (((((ab)b)c)a)a)(ab) &
  -2 (((((ab)c)a)b)a)(ab) \\[-1pt]
  -  ((((a(ab))a)c)a)(bb) &
  +2 ((((a(ab))b)a)b)(ac) &
  -2 ((((a(ab))b)b)a)(ac) \\[-1pt]
  +2 ((((a(ab))b)c)a)(ab) &
  +  ((((a(ac))a)b)a)(bb) &
  -  ((((a(ac))b)a)a)(bb) \\[-1pt]
  -2 ((((aa)(ab))a)b)(bc) &
  +2 ((((aa)(ab))b)a)(bc) &
  \end{array}
  \]
  \end{center}
  \caption{Special identity for partition 431 (terms 151 to 296)}
  \label{specialidentitypart2}
  \end{table}

Since the rank has increased by 1 for partition 431, we expect there to be a
new identity in which every monomial consists of a right-commutative
association type applied to a permutation of $aaaabbbc$. There are 106
association types and 280 permutations, so the number of monomials is at most
29680; the other partitions with new identities give larger upper bounds.
Right-commutativity implies that many of these monomials are equal; if we count
only those which equal their own straightened forms then we obtain 12131
distinct monomials. For dialgebra monomials we have 8 association types giving
2240 distinct monomials. The expansion of each right-commutative monomial is a
linear combination of 128 dialgebra monomials.

In step 1, we create a $12411 \times 12131$ matrix with a $12131 \times 12131$
upper block and a $280 \times 12131$ lower block.  For each lifted identity in
degree 8, we apply all 280 substitutions of $aaaabbbc$ for the variables, store
these nonlinear identities in the lower block, and compute the RCF using
arithmetic modulo $p = 101$. After this process is complete, the rank of the
matrix is 11020.

In step 2, we create a $2240 \times 12131$ matrix, initialize the columns with
the coefficients of the expansions of the nonlinear right-commutative
monomials, and compute the RCF using arithmetic modulo $p = 101$. The rank of
the matrix is 1110 and so the nullspace has dimension 11021.

The nullity from step 2 is exactly one more than the rank from step 2, as
expected from row 9 of Table \ref{degree8ranks}. The row space from step 1 is a
subspace of the nullspace from step 2. We need to find a nullspace vector which
is not in the row space.

In step 3, we compute the canonical basis of the nullspace from step 2. We sort
the basis vectors by increasing size; we define the ``size'' of a vector over a
finite field to be the number of distinct coefficients. We include the basis
vectors one at a time as a new bottom row of the matrix from step 1 until we
find the first basis vector that increases the rank. We multiply this basis
vector by 2 and reduce the coefficients modulo 101 using symmetric
representatives so that all the coefficients become small integers. We obtain
the 296-term identity in Tables \ref{specialidentitypart1} and
\ref{specialidentitypart2}. We expand this identity using rational arithmetic
and verify that it collapses to zero in the free associative dialgebra; this
verifies that it is a special identity over $\mathbb{Q}$.

Our new special identity involves three variables and is linear in one
variable. Therefore the obvious generalization of Macdonald's theorem
\cite{Macdonald} to quasi-Jordan algebras is not true, since our new identity
is satisfied by all special quasi-Jordan algebras but not by all quasi-Jordan
algebras. If we assume commutativity and collect terms, we obtain a polynomial
with 191 terms in the free commutative nonassociative algebra. This commutative
identity involves three variables and is linear in one variable; it is
satisfied by all special Jordan algebras since every special Jordan algebra is
a special quasi-Jordan algebra (corresponding to an associative dialgebra in
which the two operations coincide). Therefore, by Macdonald's theorem, this
commutative identity is satisfied by all Jordan algebras, and hence must be
satisfied by the Albert algebra $H_3(\mathbb{C})$.

It follows that we cannot use the Albert algebra to give a direct proof that
our new identity is not satisfied by all semispecial quasi-Jordan algebras.


\section{Conclusion}

Semispecial quasi-Jordan algebras are a natural generalization of Jordan
algebras to a noncommutative setting.  An important open problem is to
generalize classical results on free (special) Jordan algebras to semispecial
quasi-Jordan algebras. For example:
  \begin{enumerate}
  \item[$i$)] the criterion of Cohn \cite{Cohn} for a quotient of a free
      special Jordan algebra to be special, which implies that special
      Jordan algebras do not form a variety;
  \item[$ii$)] the characterization by Cohn \cite{Cohn} of the free special
      Jordan algebra on $\le 3$ generators as the symmetric elements in a
      free associative algebra;
  \item[$iii$)] the theorem of Macdonald \cite{Macdonald} on special Jordan
      identities in 3 variables;
  \item[$iv$)] the theorem of Shirshov \cite{Shirshov2} that the free
      Jordan algebra on two generators is special (see also Jacobson and
      Paige \cite{JacobsonPaige}).
  \end{enumerate}


\section*{Acknowledgements}

We thank the anonymous referee for very helpful remarks: the algorithm in
subsection \ref{referee'salgorithm} is based on the referee's report. We thank
Ra\'ul Felipe and Hader Elgendy for bringing to our attention the papers of
Kolesnikov \cite{Kolesnikov} and Bergdolt \cite{Bergdolt}. Murray Bremner
thanks NSERC for financial support through a Discovery Grant, and IME-USP for
its hospitality in May and June 2009 when most of this work was completed.



\begin{thebibliography}{99}

\bibitem{AlbertPaige}
  \textsc{A. A. Albert and L. J. Paige}:
  On a homomorphism property of certain Jordan algebras.
  \emph{Trans. Amer. Math. Soc.}
  93 (1959) 20--29.
  MR0108524 (21 \#7240)

\bibitem{Bergdolt}
  \textsc{G. Bergdolt}:
  Tilted irreducible representations of the permutation group.
  \emph{Comput. Phys. Comm.}
  86 (1995), no. 1-2, 97--104.
  MR1327568 (96e:20015)

\bibitem{Bremner}
  \textsc{M. R. Bremner}:
  On the definition of quasi-Jordan algebra.
  \emph{Comm. Algebra}
  (to appear, accepted 23 October 2009).

\bibitem{BremnerPeresi}
  \textsc{M. R. Bremner and L. A. Peresi}:
  Nonhomogeneous subalgebras of Lie and special Jordan superalgebras.
  \emph{J. Algebra}
  322, 6 (2009) 2000-2026.
  MR2542829

\bibitem{Clifton}
  \textsc{J. M. Clifton}:
  A simplification of the computation of the natural representation of
  the symmetric group $S_n$.
  \emph{Proc. Amer. Math. Soc.}
  83 (1981), no. 2, 248--250.
  MR0624907 (82j:20024)

\bibitem{Cohn}
  \textsc{P. M. Cohn}:
  On homomorphic images of special Jordan algebras.
  \emph{Canadian J. Math.}
  6 (1954) 253--264.
  MR0060496 (15,678c)

\bibitem{Glennie1}
  \textsc{C. M. Glennie}:
  \emph{Identities in Jordan Algebras}.
  Ph.D. dissertation, Yale University, 1963.

\bibitem{Glennie2}
  \textsc{C. M. Glennie}:
  Some identities valid in special Jordan algebras but not valid in all Jordan algebras.
  \emph{Pacific J. Math.}
  16 (1966) 47--59.
  MR0186708 (32 \#4166)

\bibitem{Glennie3}
  \textsc{C. M. Glennie}:
  Identities in Jordan algebras.
  \emph{Computational Problems in Abstract Algebra}.
  Proc. Conf., Oxford, 1967, Pergamon, Oxford, 1970, pages 307--313.
  MR0255629 (41 \#289)

\bibitem{Hentzel1}
  \textsc{I. R. Hentzel}:
  Processing identities by group representation.
  \emph{Computers in Nonassociative Rings and Algebras}.
  Special session, 82nd Annual Meeting Amer. Math. Soc.,
  San Antonio, 1976, pp. 13--40.
  Academic Press, New York, 1977.
  MR0463251 (57 \#3204)

\bibitem{Hentzel2}
  \textsc{I. R. Hentzel}:
  Applying group representation to nonassociative algebras.
  \emph{Ring Theory}.
  Proc. Conf., Ohio Univ., Athens, Ohio, 1976, pp. 133--141.
  Lecture Notes in Pure and Appl. Math., Vol. 25, Dekker, New York, 1977.
  MR0435159 (55 \#8120)

\bibitem{Hentzel3}
  \textsc{I. R. Hentzel}:
  Special Jordan identities.
  \emph{Comm. Algebra}
  7 (1979) no. 16, 1759--1793.
  MR0546197 (81a:17008)

\bibitem{Jacobs}
  \textsc{D. P. Jacobs}:
  Albert: an Interactive Program to Assist the Specialist in the Study of
  Nonassociative Algebra.
  \texttt{http://www.cs.clemson.edu/{\textasciitilde}dpj/albertstuff/albert.html}

\bibitem{JacobsonPaige}
  \textsc{N. Jacobson and L. J. Paige}:
  On Jordan algebras with two generators.
  \emph{J. Math. Mech.}
  6 (1957) 895--906.
  MR0092780 (19,1157b)

\bibitem{JamesKerber}
  \textsc{G. James and A. Kerber}:
  \emph{The Representation Theory of the Symmetric Group}.
  Encyclopedia of Mathematics and its Applications, 16.
  Addison-Wesley Publishing Co., Reading, Mass., 1981.
  MR0644144 (83k:20003)

\bibitem{Kernighan}
  \textsc{B. W. Kernighan and D. Ritchie}:
  \emph{The C Programming Language}.
  Second edition.
  Prentice-Hall, 1988.

\bibitem{Kolesnikov}
  \textsc{P. S. Kolesnikov}:
  Varieties of dialgebras and conformal algebras.
  \emph{Sibirsk. Mat. Zh.} 49 (2008), no. 2, 322--339;
  translation in \emph{Sib. Math. J.} 49 (2008), no. 2, 257--272.
  MR2419658 (2009b:17002)

\bibitem{Kurosh}
  \textsc{A. G. Kurosh}:
  Nonassociative free sums of algebras.
  \emph{Mat. Sb. N.S.}
  37(79) (1955) 251--264.
  MR0078356 (17,1180i)

\bibitem{Loday1}
  \textsc{J.-L. Loday}:
  Une version non commutative des alg\`ebres de Lie: les alg\`ebres de Leibniz.
  \emph{Enseign. Math.}
  (2) 39 (1993), no. 3-4, 269--293.
  MR1252069 (95a:19004)

\bibitem{Loday2}
  \textsc{J.-L. Loday}:
  Alg\`ebres ayant deux op\'erations associatives (dig\`ebres).
  \emph{C. R. Acad. Sci. Paris S\'er. I Math.}
  321 (1995), no. 2, 141--146.
  MR1345436 (96f:16013)

\bibitem{Loday3}
  \textsc{J.-L. Loday}:
  Dialgebras.
  \emph{Dialgebras and Related Operads}, 7--66,
  Lecture Notes in Math., 1763, Springer, Berlin, 2001.
  MR1860994 (2002i:17004)

\bibitem{Macdonald}
  \textsc{I. G. Macdonald}:
  Jordan algebras with three generators.
  \emph{Proc. London Math. Soc.}
  (3) 10 (1960) 395--408.
  MR0126473 (23 \#A3769)

\bibitem{Malcev}
  \textsc{A. I. Malcev}:
  On algebras defined by identities.
  \emph{Mat. Sbornik N.S.}
  26(68) (1950) 19--33.
  MR0033280 (11,414d)

\bibitem{Maple}
  \textsc{Maple}:
  Maple 13 -- The Essential Tool for Mathematics and Modeling.
  \emph{Maplesoft, a Division of Waterloo Maple Inc.}
  \texttt{http://www.maplesoft.com}

\bibitem{Pozhidaev}
  \textsc{A. P. Pozhidaev}:
  0-dialgebras with bar unity, Rota-Baxter and 3-Leibniz algebras.
  Pages 245--256 of
  \emph{Groups, Rings and Group Rings},
  edited by A. Giambruno, C. Polcino Milies and S. K. Sehgal,
  Contemporary Mathematics 499, American Mathematical Society, 2009.

\bibitem{Shirshov1}
  \textsc{A. I. Shirshov}:
  Subalgebras of free commutative and free anticommutative algebras.
  \emph{Mat. Sbornik N.S.}
  34(76) (1954) 81--88.
  MR0062112 (15,929b)

\bibitem{Shirshov2}
  \textsc{A. I. Shirshov}:
  On special $J$-rings.
  \emph{Mat. Sb. N.S.}
  38(80) (1956) 149--166.
  MR0075936 (17,822e)

\bibitem{Shirshov3}
  \textsc{A. I. Shirshov}:
  Some problems in the theory of rings that are nearly associative.
  \emph{Uspekhi Mat. Nauk}
  13 (1958), no. 6(84), 3--20.
  MR0102532 (21 \#1323);
  translation in
  \emph{Non-associative Algebra and its Applications}, 441--459,
  Lect. Notes Pure Appl. Math., 246,
  Chapman \& Hall/CRC, Boca Raton, 2006.

\bibitem{Sloane}
  \textsc{N. J. A. Sloane}:
  \emph{The On-Line Encyclopedia of Integer Sequences}.
  AT\&T Labs - Research.
  \texttt{http://www.research.att.com/{\textasciitilde}njas/sequences}

\bibitem{Specht}
  \textsc{W. Specht}:
  Gesetze in Ringen, I.
  \emph{Math. Z.}
  52 (1950) 557--589.
  MR0035274 (11,711i)

\bibitem{VelasquezFelipe1}
  \textsc{R. Vel\'asquez and R. Felipe}:
  Quasi-Jordan algebras.
  \emph{Comm. Algebra}
  36 (2008), no. 4, 1580--1602.
  MR2410352 (2009b:17004)

\bibitem{Zhevlakov}
  \textsc{K. A. Zhevlakov, A. M. Slinko, I. P. Shestakov and A. I. Shirshov}:
  \emph{Rings That Are Nearly Associative}.
  Translated from the Russian by Harry F. Smith.
  Pure and Applied Mathematics, 104.
  Academic Press, Inc., New York, 1982.
  MR0668355 (83i:17001)

\end{thebibliography}
\end{document}